\documentclass{IEEEtran}
\usepackage{cite}
\usepackage{amsmath,amssymb,amsfonts}
\usepackage{algorithmic}
\usepackage{graphicx}
\usepackage{textcomp}
\def\BibTeX{{\rm B\kern-.05em{\sc i\kern-.025em b}\kern-.08em
    T\kern-.1667em\lower.7ex\hbox{E}\kern-.125emX}}

\usepackage{subfig}

\begin{document}
\title{Fully Automated Adaptive Parameter Selection for 3-D High-order Nyström Boundary Integral Equation Methods}
\author{Davit Aslanyan, \IEEEmembership{Student Member, IEEE}, and Constantine Sideris, \IEEEmembership{Senior Member, IEEE}

\thanks{The authors gratefully acknowledge support by the Air Force Office of Scientific Research (FA9550-20-1-0087 and FA9550-25-1-0020) and the National Science Foundation (CCF-2047433).}

\thanks{D. Aslanyan is with the Department of Electrical and Computer Engineering, University of Southern California, Los Angeles, CA 90089, USA and C. Sideris is with the Department of Electrical Engineering, Stanford University, Stanford, CA 94305 (e-mails: daslanya@usc.edu, sideris@stanford.edu).}}

\maketitle

\begin{abstract}
We present an adaptive Chebyshev-based Boundary Integral Equation (CBIE) solver for electromagnetic scattering from smooth perfect electric conductor (PEC) objects. The proposed approach eliminates manual parameter tuning by introducing (i) a unified adaptive quadrature strategy for automatic selection of the near-singular interaction distance and (ii) an adaptive computation of all self- and near-singular precomputation integrals to a prescribed accuracy using Gauss–Kronrod (h-adaptive) or Clenshaw–Curtis (p-adaptive) rules and singularity-resolving changes of variables. Both h-adaptive and p-adaptive schemes are explored within this framework, ensuring high-order accuracy and robustness across a broad range of geometries without loss of efficiency. Numerical results for canonical and complex CAD geometries demonstrate that the adaptive solver achieves accuracy and convergence rates comparable to optimally tuned fixed-grid CBIE implementations, while offering automation and scalability to electrically large, geometrically complex problems.
\end{abstract}


\section{Introduction}
\label{sec:introduction}
\IEEEPARstart{E}{fficient} and accurate electromagnetic (EM) simulations for solving Maxwell's equations are crucial for advancements in modern technology, including communication systems, medical imaging, nanophotonics, and aerospace engineering. EM scattering problems only have analytical solutions for some canonical objects, such as infinite cylinders or spheres. However, for most realistic electrically large objects of interest (e.g., aircraft, circuit interconnects, nanophotonic devices), scattering problems can only be solved numerically, highlighting the importance and necessity of efficient and accurate solvers.

The Finite Difference Method (FDM, \cite{zhou1993finite}) and the Finite Element Method (FEM, \cite{dhatt2012finite}) are two examples of numerical solvers for such problems, which have widespread use and popularity due to their relative simplicity, robustness, and ease of implementation. However, FDM and FEM are both volumetric methods and require meshing of the entire simulation domain. They also both suffer from numerical dispersion and low accuracy due to local truncation error from low-order approximations of the derivative operators, and typically require extensive geometry processing with complicated volumetric meshes to faithfully represent complex CAD models (\cite{warren_fem_numerical_dispersion_1995, warren_fem_numerical_dispersion_2002, pereda_fdtd_dispersion_1998,fdtd_review_nature, fdtd_geometry_processing_2024}). 

Boundary Integral Equation (BIE) methods constitute another popular solution method and offer a number of advantages over FDM and FEM. BIE methods require meshing only the boundaries between different domains (requiring 2D surface meshes to represent 3D objects) and, as such, can significantly reduce the number of unknowns necessary to achieve a similar accuracy when compared to volumetric approaches. BIE approaches can more readily be integrated with CAD models and support high-order curvilinear NURBS-based meshes due to the considerably simpler 2D meshing required~\cite{bie_cad_optics,jin_density_continuity,jin_ojap,simpson_bie,li_cad_bie}. They do not suffer from numerical dispersion due to analytical propagation from sources to targets using the Green's function, and, when suitably accelerated via algorthmic and hardware acceleration approaches, can achieve high efficiency in terms of computing times and parallel performance \cite{song1995mlfmm,song1997mlfmm,jagabandhu_ifgf,bruno2001mlfmm}.

The Method of Moments (MoM)~\cite{harrington1993field} is perhaps one of the most well known discretization approaches for BIEs, in which unknown surface currents are expanded in terms of basis functions and tested against weighing functions. Low-order implementations, most notably those based on Rao–Wilton–Glisson (RWG) basis functions \cite{rao1982electromagnetic}, are robust and divergence-conforming but require very fine meshes to capture geometries with curvature or electrically large surfaces, leading to high computational costs. They also only exhibit linear convergence in the solution accuracy with respect to the mesh resolution. Higher-order MoM schemes~\cite{peterson1998computational} mitigate this by using high-order basis functions, but often come with significant implementation complexity and computational costs associated with the 4D integrals required for MoM.  

Nyström-based BIE methods discretize the integral equations directly and use ``point-matching" for testing, bypassing the need for basis and testing functions. Locally corrected Nyström (LCN) methods \cite{peterson2010introduction} have been widely used in the literature to achieve high-order convergence by introducing local kernel corrections near singularities. While accurate, the construction of these correction terms can be complicated, computationally inefficient, and potentially numerically unstable for higher orders, especially for vector electromagnetic formulations on complex geometries \cite{shafieipour2014efficient_lcn,shafieipour2013equivalence_lcn}.  

The Chebyshev-based Boundary Integral Equation (CBIE) method is a highly efficient implementation of the Nyström method ~\cite{bruno2020chebyshev,hu2021chebyshev}. The CBIE method represents the geometry of the objects of interest with a set of nonoverlapping, curvilinear, quadrilateral patches, which enables an accurate representation of complex curvatures with very coarse meshes.  On each patch, the surface density unknowns are represented on a 2D Chebyshev grid using a tensor product of the nodes of 1D Chebyshev polynomials. The method leverages the fact that the unknowns are located on Chebyshev grids to efficiently calculate interactions with smooth kernels (``far" interactions) using F\'ejer's first quadrature rule. However, when target point of interest is near (``near-singular" interaction) or on (``singular" interaction) the integration patch, the integral operator kernels can exhibit near-singular or singular behavior. The CBIE method treats these interactions by: (i) applying a singularity-canceling Change-of-Variables (CoV) to calculate integration weights (integral operator kernel's action against Chebyshev polynomials) via a precomputation step on refined grids with high accuracy, (ii) expanding the unknown densities in terms of Chebyshev polynomials, and (iii) calculating the action of the integral operator's kernel against the density by multiplying the precomputed quadrature weights against the corresponding Chebyshev coefficients of the density and accumulating them. Compared to low-order RWG-based MoM, CBIE achieves higher accuracy and order of convergence with fewer unknowns; compared to LCN, CBIE can calculate the correction terms considerably more efficiently~\cite{jin_hmatrix}.  

Despite these advantages, existing CBIE implementations suffer from two serious drawbacks. First, the optimal choice of the ``near-singular" interaction distance parameter $\Delta_{\mathrm{near}}$ can strongly depend on the particular problem and choice of discretization. 
$\Delta_{\mathrm{near}}$ determines whether the method uses F\'ejer's first quadrature rule (for target points that are further away than $\Delta_{\mathrm{near}}$ from a source patch) or utilizes Chebyshev expansion of the density and the precomputation weights to calculate the integrals accurately otherwise. The $\Delta_{\mathrm{near}}$ parameter can take a wide range of values, making its selection non-intuitive and requiring extensive manual tuning for optimal performance. Second, the choice of the parameter $N_{\beta}$ (which determines the number of integration points in each $u$ and $v$ direction of the refined quadrature grid for evaluation of precomputation weights to the desired accuracy), also requires parameter tuning to achieve optimal performance. Furthermore, the optimal $\Delta_{near}$ and $N_{\beta}$ can vary significantly from patch to patch over the same geometry, although previous CBIE implementations have all used global ``worst-case" values. These limitations can hinder robustness, scalability, and often even computational efficiency for realistic geometries.  

Recently, in the context of conference proceedings, we explored both Gauss–Kronrod~\cite{notaris2016gauss} and Clenshaw–Curtis~\cite{clenshaw1960method} adaptive quadrature rules to automatically calculate the  precomputation integrals required by the CBIE method to a desired tolerance~\cite{davit_aps2023,davit_aces2023}. In this work, we extend those methodologies, introduce a new robust approach for automatically selecting $\Delta_{\mathrm{near}}$ on a per-patch basis, and propose an adaptive integration strategy for the CBIE framework that removes all manual parameter tuning. The main contributions are:
\begin{itemize}
    \item A unified adaptive quadrature approach that automatically determines the optimal ``near-singular" interaction distance $\Delta_{\mathrm{near}}$ for each integration patch.  
    \item Adaptive computation of all precomputation integrals to the prescribed accuracy using Gauss–Kronrod and Clenshaw–Curtis rules combined with singularity-resolving changes of variables.  
    \item Negligible computational overhead: The approach delivers accuracy and efficiency comparable to, or better than, existing CBIE implementations across a wide range of scattering problems.  
\end{itemize}

The remainder of the paper is organized as follows: Section ~\ref{prelim} defines the integral equation formulation used in this work and provides a brief overview of the CBIE method, Section \ref{adaptiveprecomp} introduces the two-stage \textit{h}- and \textit{p}- adaptive integration approaches, discusses the singularity-canceling Change-of-Variables utilized in the adaptive schemes, and compares the computational cost of each approach. Section \ref{results} presents numerical examples that demonstrate the robustness of the proposed approach. The Magnetic Field Integral Equation (MFIE) is solved on two PEC canonical geometries (a sphere and a toroid), as well as a complex CAD model of a glider. The results show that the adaptive strategy enables a robust, high-order error-controlled CBIE solver that is completely free of any manual parameter tuning and scalable to complex, electrically large geometries. We note that all of the techniques introduced in this work can be directly applied to any other integral formulations, including the Electric Field Integral Equation (EFIE), Combined Field Integral Equation (CFIE), as well as formulations for modeling dielectric materials, such as the PMCHWT and M\"uller formulations.

\section{Preliminaries}\label{prelim}

\subsection{Integral Equation Formulation}

\begin{figure}[!t]
\centerline{\includegraphics[width=0.5\columnwidth]{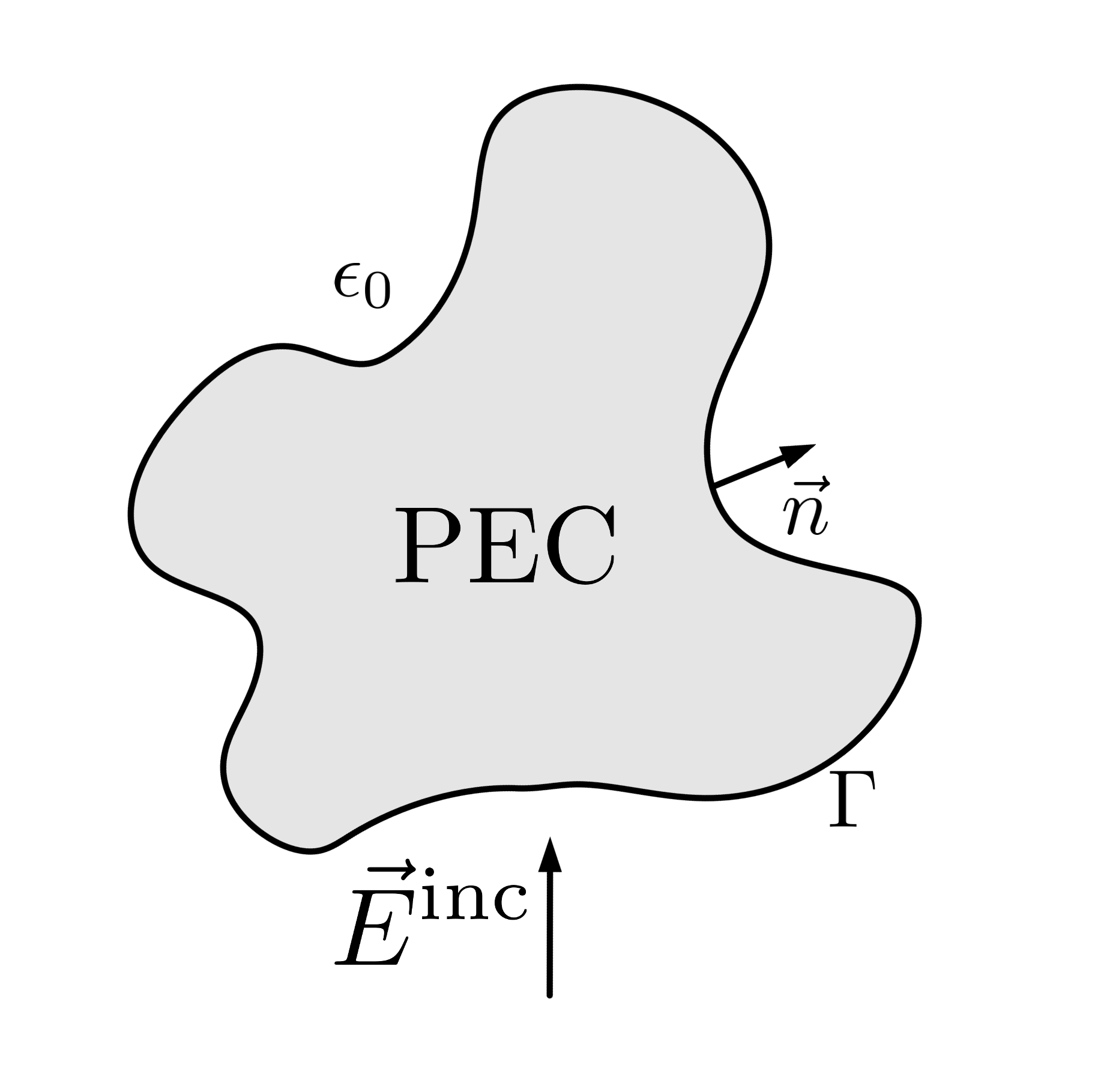}}
\caption{Plane wave incident upon a PEC object.}
\label{fig:incident_field}
\end{figure}

We consider the problem of electromagnetic scattering by closed PEC objects, enclosed by a boundary surface $\Gamma$ in 3D space. The objects are illuminated by an incident plane wave
\[
\Vec{E}^{\textrm{inc}} = e^{ik_0 z}\hat{x},
\]
with free-space wavenumber $k_0$ (Fig. \ref{fig:incident_field}).

For the sake of simplicity, we choose to use the MFIE to model this problem, although any other suitable formulation, such as the EFIE or CFIE can be readily used. The integral equation expresses the scattered fields in terms of the induced surface current density $\Vec{J} = \hat{n} \times \Vec{H}$.

The classical MFIE is expressed as
\begin{equation}
    \frac{1}{2}\Vec{J} + K[\Vec{J}] = \hat{n} \times \Vec{H}^{\textrm{inc}},
\end{equation}
where the operator $K$ is defined by
\begin{equation}
    K[\Vec{a}](\mathbf{r}) = \hat{n}(\mathbf{r}) \times \int_\Gamma \Vec{a}(\mathbf{r'}) \times \nabla G(\mathbf{r}-\mathbf{r'})\, d\sigma(\mathbf{r'}).
\end{equation}
Here $G(\mathbf{r}-\mathbf{r'}) = \tfrac{e^{-ik_0|\mathbf{r}-\mathbf{r'}|}}{4\pi|\mathbf{r}-\mathbf{r'}|}$ is the free-space scalar Green's function for the Helmholtz equation, $\mathbf{r}$ denotes the target point, $\mathbf{r'}$ the integration variable, $\nabla$ the gradient with respect to $\mathbf{r}$, and $\Gamma$ the surface of the object.

Following \cite{garza2020boundary}, the operator $K$ can be reformulated using vector calculus identities and by decomposing the gradient into tangential and normal components. This yields a representation in terms of the weakly singular single-layer and adjoint double-layer operators:
\begin{equation}\label{singlelayer}
     S[\Vec{a}](\mathbf{r}) = \int_\Gamma G(\mathbf{r}-\mathbf{r'})\Vec{a}(\mathbf{r'})\, d\sigma(\mathbf{r'}), 
\end{equation}
\begin{equation}\label{adjointdoublelayer}   
    D[\Vec{a}](\mathbf{r}) = \int_\Gamma \frac{\partial G(\mathbf{r}-\mathbf{r'})}{\partial \Vec{n}(\mathbf{r})}\Vec{a}(\mathbf{r'})\, d\sigma(\mathbf{r'}).
\end{equation}
In this notation, $K$ can be expressed as
\begin{multline}\label{eq:final_mfie_kernel}
    K[\Vec{a}](\mathbf{r}) = D[\Vec{a}](\mathbf{r}) - \\ \left( \Vec{e}^{u}(\mathbf{r})\left[ \Vec{n}(\mathbf{r})\cdot\frac{\partial}{\partial u}\right] + \Vec{e}^{v}(\mathbf{r})\left[ \Vec{n}(\mathbf{r})\cdot\frac{\partial}{\partial v}\right] \right) S[\Vec{a}](\mathbf{r}), 
    \quad \mathbf{r} \in \Gamma,
\end{multline}
where $\Vec{e}^u$ and $\Vec{e}^v$ are contravariant basis vectors, and $\Vec{n}$ is the outward-pointing unit normal. The advantage of this formulation is that when discretizing using the CBIE method, it requires only precomputing two scalar operators ($S$ and $D$) for the singular and near-singular interactions compared to the four quantities required otherwise, leading to significant savings in storage and computational effort.

\subsection{Chebyshev-based Boundary Integral Equation Framework}\label{cbie_background}

To discretize the MFIE, we employ the Chebyshev-based Nystr\"om (CBIE) method introduced in~\cite{bruno2020chebyshev,hu2021chebyshev}. This subsection briefly reviews the CBIE framework and establishes the foundation for the adaptive precomputation approach presented in Section~\ref{adaptiveprecomp}. 

\begin{figure}[!t]
\centerline{\includegraphics[width=0.8\columnwidth]{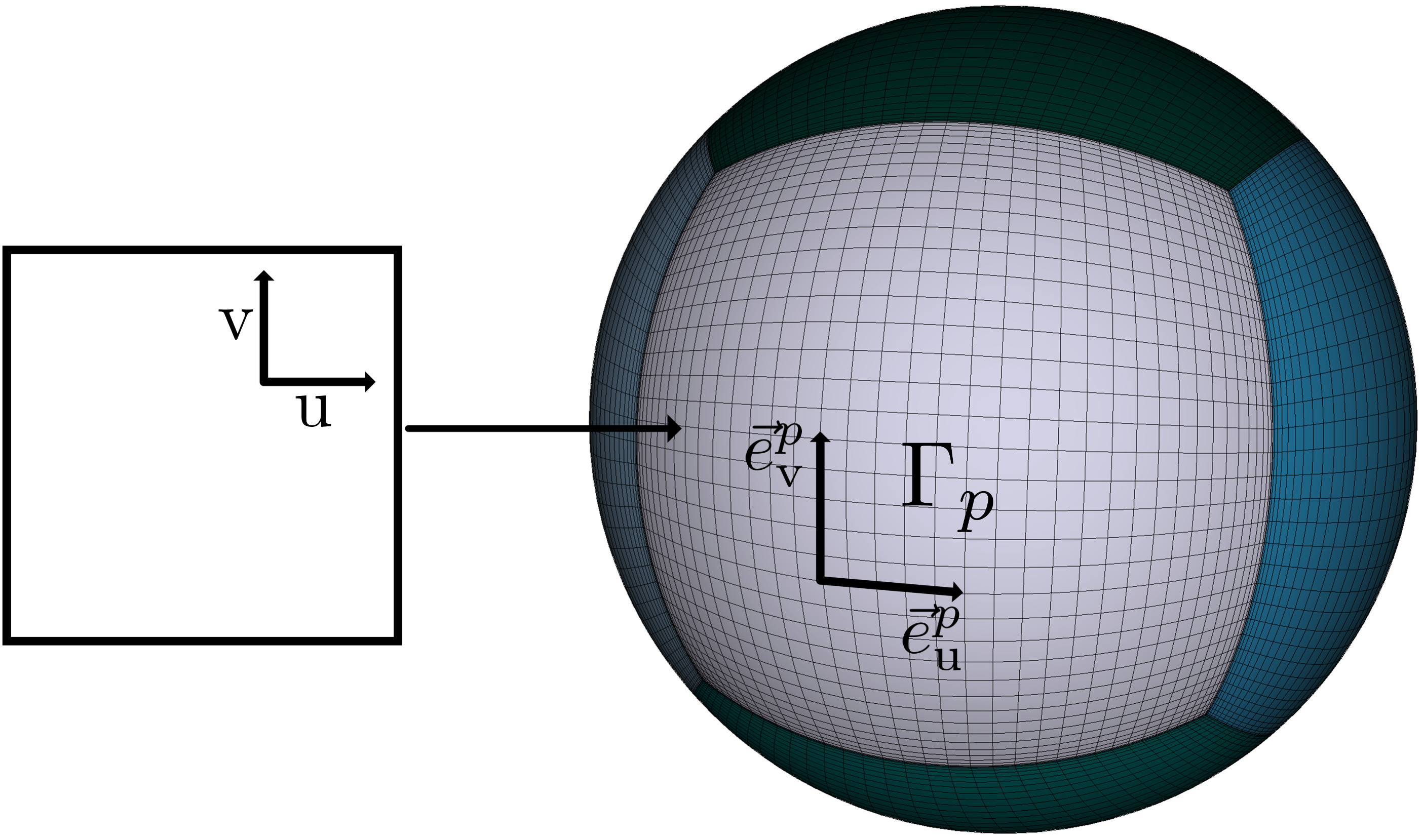}}
\caption{Parametric mapping and surface partitioning.}
\label{fig:patch_splitting}
\end{figure}

The CBIE method partitions the surface $\Gamma$ into $\mathrm{M}$ non-overlapping quadrilateral patches (Fig \ref{fig:patch_splitting}):
\begin{equation}\label{partition}
\Gamma = \bigcup_{p=1}^M \Gamma^p,
\end{equation}
where each patch $\Gamma^p$ is parametrized by a mapping
\begin{equation}\label{param}
\mathbf{r}_p(u,v):[-1,1]\times[-1,1] \to \Gamma^p \subset \mathbb{R}^3.
\end{equation}

For each patch $p$, the differential geometric quantities are defined as follows. The tangential basis vectors are
\[
\Vec{e}_u^p = \frac{\partial \mathbf{r}_p(u,v)}{\partial u}, 
\quad
\Vec{e}_v^p = \frac{\partial \mathbf{r}_p(u,v)}{\partial v},
\]
the outward-pointing unit normal is
\[
\Vec{n}_p(\mathbf{r}) = \frac{\Vec{e}_u^p \times \Vec{e}_v^p}{\|\Vec{e}_u^p \times \Vec{e}_v^p\|},
\]
and the metric tensor is
\begin{equation}\label{metric_tensor}
G^p = \begin{bmatrix}
g^p_{uu} & g^p_{uv} \\
g^p_{vu} & g^p_{vv}
\end{bmatrix}, \quad 
g_{ij}^p = \Vec{e}_i^p \cdot \Vec{e}_j^p, \; i,j \in \{u,v\}.
\end{equation}
The surface element can then be defined as $d\sigma^p = \sqrt{|G^p|}\,du\,dv$, where $|G^p|$ denotes the determinant of the metric tensor. The contravariant basis vectors are given by $\Vec{e}^i_p = g^{ij}_p \Vec{e}_j^p$, with $g^{ij}_p$ the components of the inverse of $G^p$. The surface current density can be expanded as
\[
\Vec{J} = J^u \Vec{e}_u + J^v \Vec{e}_v.
\]

\begin{figure*}[!ht]
    \centering
    \subfloat[]{\label{fig:cov_single}\includegraphics[width=0.33\textwidth]{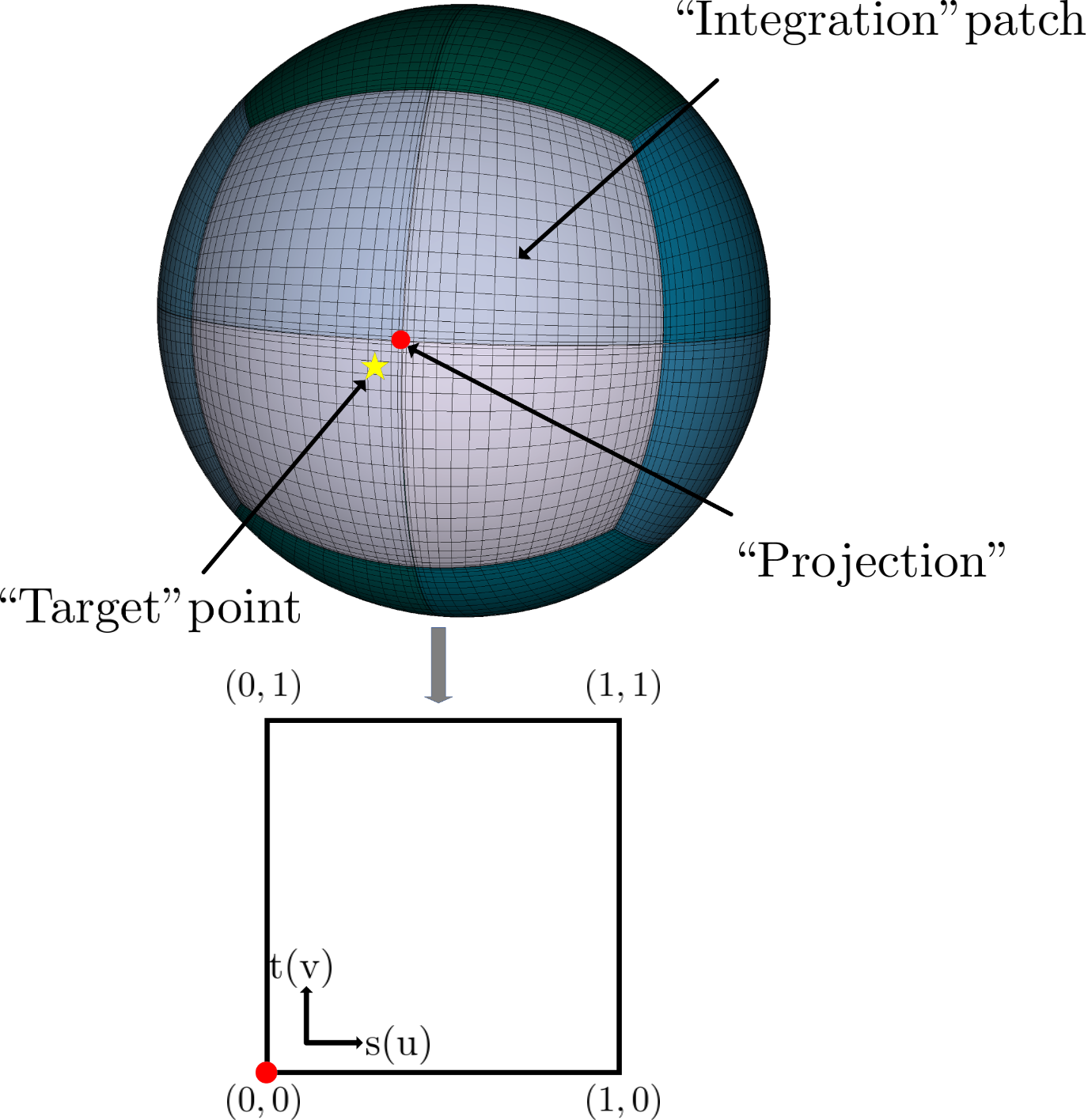}}
    \subfloat[]{\label{fig:cov_two}\includegraphics[width=0.33\textwidth]{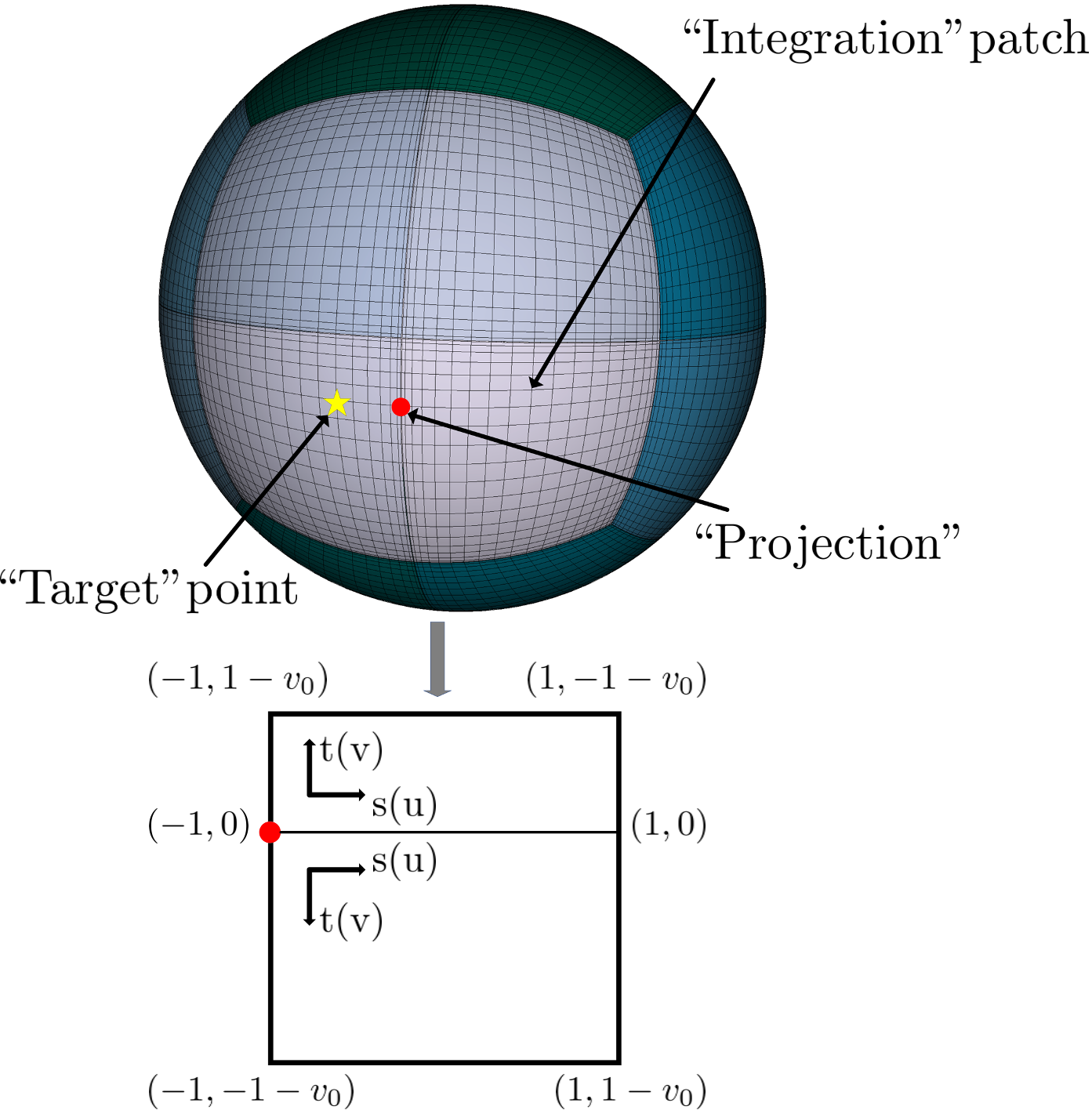}}
    \subfloat[]{\label{fig:cov_four}\includegraphics[width=0.33\textwidth]{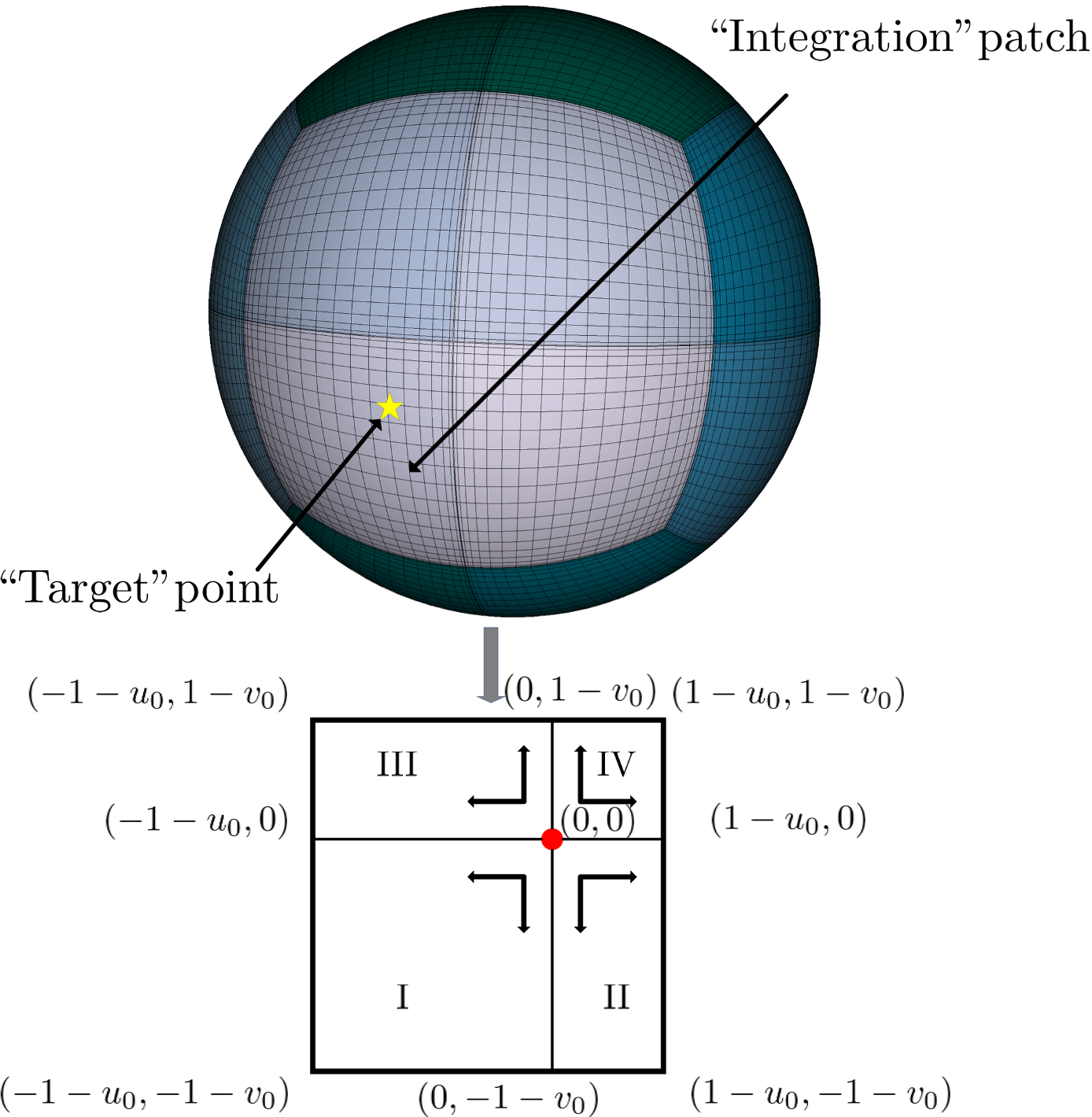}}
    
    \caption{Singular Change-of-Variables depending on the location of the ``target" point and the ``integration" patch. (a) Integration patch not subdivided: CoV applied on a single subpatch, (b) Integration patch divided in $v$-direction: CoV applied on two subpatches, (c) Integration patch divided in both $u$ and $v$ directions: CoV applied on four subpatches. Note that in this case, the ``projection" point coincides with the ``target" point.}
    \label{fig:cov}
\end{figure*}

Any function $f(\mathbf{r})$ defined on $\Gamma^p$ corresponds to $f(\mathbf{r}_p(u,v))$ based on the patch parametrization, and is discretized on a tensor-product Chebyshev grid~\cite{chebyshevnodes} given by:
\begin{multline}\label{negoneonenodes}
x_j = \cos \left( \pi \frac{2j + 1}{2Q} \right), \quad j = 0, \ldots, Q - 1, \\[4pt]
w_j = \frac{2}{Q} \left( 1 - 2 \sum_{k=1}^{Q/2} \frac{1}{4k^2 - 1} \cos \left( k\pi \frac{2j + 1}{Q} \right) \right),
\end{multline}
where $Q \in \{n,m\}$ are the numbers of discretization points in the $u$ and $v$ directions. While it is not required by the CBIE framework, we choose $n=m$ for simplicity. 

The action of an integral operator with kernel $B(\mathbf{r},\mathbf{r}')$ on $f$ can be simplified as a summation of patchwise computed integrals:
\begin{equation}
    I(\mathbf{r}) = \int_{\Gamma} B(\mathbf{r},\mathbf{r}')f(\mathbf{r}')\, d\sigma(\mathbf{r}')
    = \sum_{p=1}^M I^p(\mathbf{r}),
\end{equation}
where
\begin{equation}
 I^p(\mathbf{r}) = \int_{-1}^1\!\!\int_{-1}^1 B(\mathbf{r},\mathbf{r}_p(u,v)) f(\mathbf{r}_p(u,v)) \sqrt{|G^p|}\,du\,dv.
\end{equation}

Depending on the relative location of $\mathbf{r}$ and $\Gamma^p$, the kernel $B$ may be smooth or singular. The CBIE framework achieves high-order accuracy by treating these cases separately.  

\textit{Far interactions}: If $\mathbf{r}$ is sufficiently far from $\Gamma^p$, the integral is approximated using F\'ejer's first quadrature rule:
\begin{equation}\label{eq:farinteractions}
 I^p(\mathbf{r}) \approx \sum_{i=0}^{n-1}\sum_{j=0}^{m-1} B(\mathbf{r},u_i,v_j)f(u_i,v_j)\sqrt{|G^p(u_i,v_j)|}\, w_i w_j.
\end{equation}

\textit{Near interactions}: If $\mathbf{r}$ lies on or near $\Gamma^p$, the function $f$ is expanded in Chebyshev polynomials:
\begin{equation}\label{eq:chebyexpansion}
f(u,v) \approx \sum_{i=0}^{n-1}\sum_{j=0}^{m-1} a_{i,j}^p T_i(u)T_j(v),
\end{equation}
with coefficients
\begin{equation}\label{eq:cheby_expansion_coefficients}
a_{i,j}^p = \frac{\alpha_i\alpha_j}{nm} \sum_{l=0}^{n-1}\sum_{t=0}^{m-1} f(u_l,v_t)T_i(u_l)T_j(v_t),
\end{equation}
where $\alpha_i = 1$ if $i=0$ and $\alpha_i=2$ otherwise.  
The operator action is then computed as
\begin{equation}\label{eq:nearinteraction}
I^p(\mathbf{r}) \approx \sum_{i=0}^{n-1}\sum_{j=0}^{m-1} a_{i,j}^p \beta_{i,j}^p,
\end{equation}
where $\beta_{i,j}^p$ are the precomputation weights:
\begin{equation}\label{eq:nearinteraction_weights}
 \beta_{i,j}^p = \int_{-1}^1\!\!\int_{-1}^1 B(\mathbf{r},\mathbf{r}_p(u,v))T_i(u)T_j(v)\sqrt{|G^p|}\,du\,dv.
\end{equation}

In this work, the MFIE is implemented using only the adjoint double layer operator, the single-layer operator, and their numerical derivatives. This minimizes the number of precomputations and the kernel $B$ reduces to
\begin{equation}\label{kerneldef}
 B(\mathbf{r},\mathbf{r}') = 
   \begin{cases}
   G(\mathbf{r}-\mathbf{r'}), & \text{for Eq.~\eqref{singlelayer}}, \\[4pt]
   \tfrac{\partial G(\mathbf{r}-\mathbf{r'})}{\partial \Vec{n}(\mathbf{r})}, & \text{for Eq.~\eqref{adjointdoublelayer}}.
   \end{cases}
\end{equation}

\section{Adaptive Integration Approach}\label{adaptiveprecomp}

An adaptive quadrature is used to achieve the desired accuracy for the precomputation of integral equation kernels against Chebyshev polynomials for the singular and near singular interactions. Target points considered potentially near a source patch are ordered from nearest to furthest, and precomputation integrals are computed for each one and compared against using the F\'ejer quadrature for far interactions. Once the error of the F\'ejer quadrature compared to the more accurate near quadrature method drops below the desired tolerance, the remainder of further away target points are considered ``far" interactions and precomputations are not computed for them. The error in both the calculation of the precomputation integrals and the determination of the near vs. far target points is estimated by calculating the action of the kernel against an ``auxiliary" density function. Since these quantities must be determined before solving the problem, the true solution density is not available. We propose a constructed ``auxiliary" density based on a plane-wave with slightly larger wavenumber than what is being solved, which helps control the accuracy of the precomputation stage and establishes the near/far interaction boundary in a consistent fashion. This function is chosen as a reliable heuristic to mimic the expected behavior of physical surface densities and thus provides a consistent basis for error estimation.

\subsection{Adaptive Framework}
The algorithm proceeds as follows:
\begin{itemize}
    \item {Precomputation accuracy:} For each integration source patch, the target points are ordered by distance. For each target point in sequence, the precomputation weights (Eq.~\eqref{eq:nearinteraction_weights}) for the action of the integral operator kernel against Chebyshev polynomials are calculated using either $h$-adaptive Gauss-Kronrod quadrature or $p$-adaptive Clenshaw-Curtis quadrature. Both adaptive quadrature schemes iteratively refine the quadrature grids until the desired accuracy tolerances are met. Sections \ref{GKsection} and \ref{CCsection} discuss the details of the implementation of the schemes.
    \item {Near/far interaction cutoff determination:} Once precomputations for each target point--integration patch interactions are calculated to the desired accuracy, the algorithm determines the ``near" interaction cutoff distance by checking the calculation error. The error is defined as the relative error between the integration results \eqref{eq:farinteractions} and \eqref{eq:nearinteraction} for the ``auxiliary" density. If the error is below the specified tolerance, all subsequent target points (which due to the sorting mentioned above are located farther away from the integration patch) are treated as ``far" interactions, and the ``near" interaction tolerance is determined.
    
\end{itemize}

\subsection{Singularity-Resolving Change of Variables}\label{sub:cov_disc}

To effectively treat operator singularities that arise when the ``target" point is on or near the integration patch, we employ a regularizing $p$-th order change of variables (Eq.\ref{eq:cov}) centered at the point ($u_0,v_0$): 
\begin{equation}\label{eq:cov}
s = \eta_s^q(u), \quad t = \eta_t^q(v),
\end{equation}
where $\eta^q$ is chosen such that its first $p-1$ derivatives vanish at the origin. This clustering of quadrature nodes around the singularity point regularizes the integrand and accelerates convergence.

The point ($u_0,v_0$) is determined by the relative position of the ``target" point and the integration patch. When the observation point is on the integration patch, ($u_0,v_0$) is set to the observation point itself, since $|r-r'|\to0 \text{\ as \ }r\to r'$. When the ``target" point is not on the integration patch but is close to it such that the kernels exhibit near-singular behavior, ($u_0,v_0$) corresponds to the point on the integration patch that minimizes the Euclidean distance to the ``target" point. 

\begin{figure*}[!t]
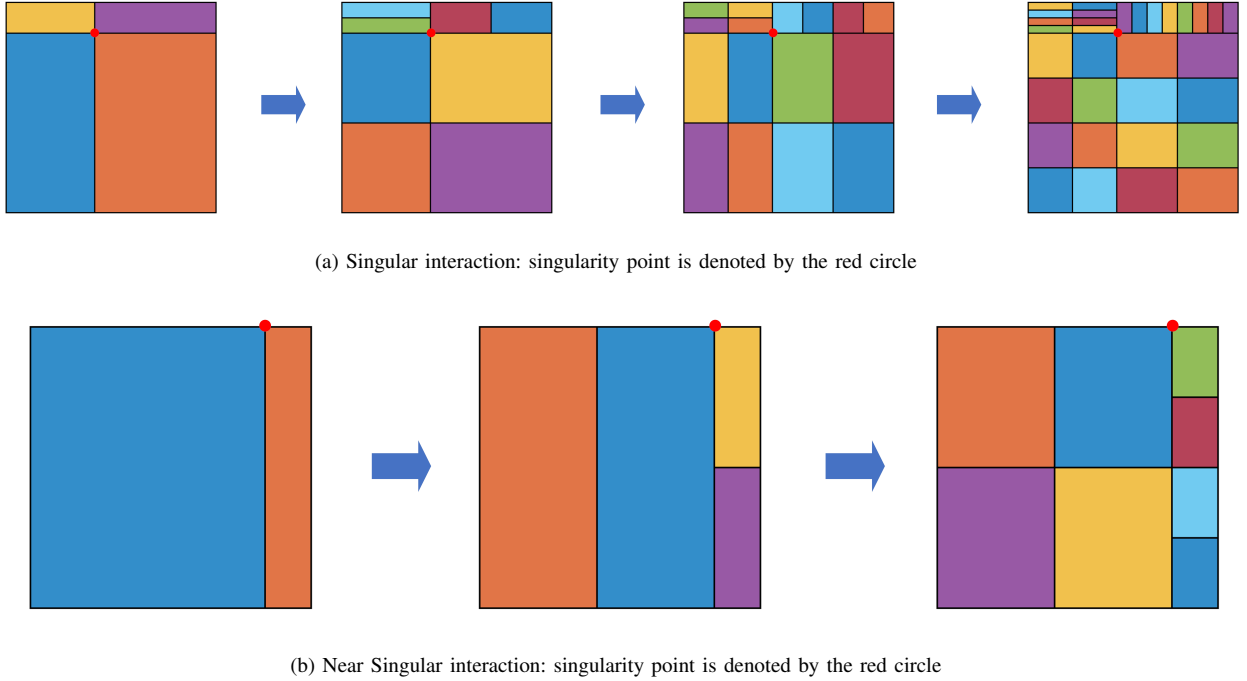

    \centering
    \subfloat[Singular interaction: singularity point is denoted by the red circle]{\label{fig:gk_splitting}\includegraphics[width=1.0\textwidth]{Figures/gk_intervals.png}}
    
    \subfloat[Near Singular interaction: singularity point is denoted by the red circle]{\label{fig:gk_splitting_near}\includegraphics[width=1.0\textwidth]{Figures/gk_intervals_near.png}}

    \caption{An example of the interval splitting used by the Gauss-Kronrod Quadrature for calculating precomputation integrals}
\end{figure*}

Each integration patch is subdivided into up to four subpatches, with ($u_0,v_0$) mapped to the origin of each subpatch. If the integration patch is divided into four subpatches the CoV (Fig. \ref{fig:cov}) has the form: 
\begin{equation}
s = 
\begin{cases}
    u^p (-u_0-1) + u_0 & \text{Interval I,}\\
    u^p (1-u_0) + u_0 & \text{Interval II,}\\
    u^p (-u_0-1) + u_0 & \text{Interval III,}\\
    u^p (1-u_0) + u_0 & \text{Interval IV,}
\end{cases}
\end{equation}
and
\begin{equation}
t = 
\begin{cases}
    v^p (-v_0-1) + v_0 & \text{Interval I,}\\
    v^p (-v_0-1) + v_0 & \text{Interval II,}\\
    v^p (1-v_0) + v_0 & \text{Interval III,}\\
    v^p (1-v_0) + v_0 & \text{Interval IV.}
\end{cases}
\end{equation}
Depending on the distance between the coordinates of the point ($u_0$, $v_0$) and domain endpoints $\{-1,1\}$, the adaptive integration approach decides whether to subdivide the integration patch into one, two, or four subpatches:
\begin{itemize}
    \item If both projection coordinates are (within $10^{-10}$ tolerance) close to either endpoint, the patch is not divided (Fig. ~\ref{fig:cov_single}),
    \item If only one of projection coordinates is (within $10^{-10}$ tolerance) close to either endpoint, the patch is divided into two subpanels (Fig. ~\ref{fig:cov_two}),
    \item If neither projection coordinate is (within $10^{-10}$ tolerance) close to either endpoint, the patch is divided into four subpanels (Fig. ~\ref{fig:cov_four}).
\end{itemize}

Although Gauss-Kronrod quadrature can integrate weak endpoint singularities without modification, incorporating this CoV improves its efficiency by smoothing out the integrand. For Clenshaw-Curtis quadrature, the CoV is indispensable, since the polar-rectangular integration converges considerably faster once the kernel singularity is smoothed~\cite{bruno2020chebyshev}. Importantly, both quadrature methods are always applied on the same CoV-transformed subpatches, ensuring consistency across the two schemes. 

\begin{figure}[!t]
\centerline{\includegraphics[width=0.75\columnwidth]{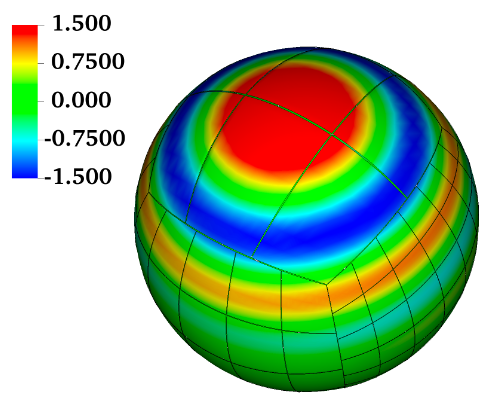}}
\caption{The real part of $J_x$ solved surface current density on a PEC Sphere of $D = 4\lambda$ diameter.}

\label{fig:sphere_density}
\end{figure}

\subsection{Gauss-Kronrod Quadrature}\label{GKsection}
The $h$-adaptive strategy is based on the nested Gauss-Kronrod (GK) rule \cite{notaris2016gauss}. In one dimension, the lower order Gauss (with $N_G$ points) and higher order Kronrod (with $N_K$ points) rules share nodes, with Gauss nodes corresponding to every-other-node of the Kronrod rule. This enables simultaneous numerical evaluation of the integral using both the Gauss and Kronrod quadrature rules using the same kernel evaluations at the Kronrod nodes. Their difference provides an error estimate at no significant additional cost due to the higher order of the Kronrod rule. This can be extended to two-dimensional integration by using the tensor-product of two 1D rules. Depending on whether Gauss (G) or Kronrod (K) rules are used in each direction we can obtain four estimates for the integrals in $(u,v)$: $
I_{GG}, \ I_{GK}, \ I_{KG}, \ I_{KK},
$
where the subscripts denote whether Gauss (G) or Kronrod (K) rules are used in each direction~\cite{davit_aps2023}. The $I_{KK}$ estimate is used as the reference as it is expected to be the most accurate. The algorithm iteratively subdivides each integration panel and applies the same order Gauss-Kronrod rule on each resulting subpanel until global absolute and/or relative error tolerances are met or the maximum number of integration panels is reached:
\[
|I_{KK} - I_{GG}| < \max(\text{tol}_{\text{abs}}, \text{tol}_{\text{rel}} \cdot |I_{KK}|).
\]

Traditional implementations of the GK quadrature uniformly divide each integration panel into 4 sub-panels, however this approach does not use all of the information available (namely the estimates $I_{GK}$ and $I_{KG}$). Our implementation considers the error estimates $|I_{KK} - I_{GK}|$ and $|I_{KK} - I_{KG}|$ to determine which direction has worse accuracy and preferentially splits in that direction, leading to an anisotropic adaptive integration approach. This leads to a more efficient allocation of computational effort since the isotropic approach always splits in both directions even if the prescribed error tolerance has been achieved in one of them.

The integrals $I_{XY}$ are defined as
\[
I_{XY}(\mathbf{r}) \approx \sum_{i=0}^{n-1}\sum_{j=0}^{m-1} a_{i,j}\beta_{i,j}^{XY},
\]
where $X,Y \in \{G,K\}$, $a_{i,j}$ are the Chebyshev expansion coefficients for the ``auxiliary" density, and $\beta_{i,j}^{XY}$ denotes the result of the precomputation integral ~\eqref{eq:nearinteraction_weights}. An example workflow of the implementation for self-singular and near-singular interactions is presented in Fig. ~\ref{fig:gk_splitting} and Fig. ~\ref{fig:gk_splitting_near}.

\begin{figure}[!ht]
    \centering
    \subfloat[]{\label{fig:sphere_rcs_vs_unknowns}\includegraphics[width=1.0\columnwidth]{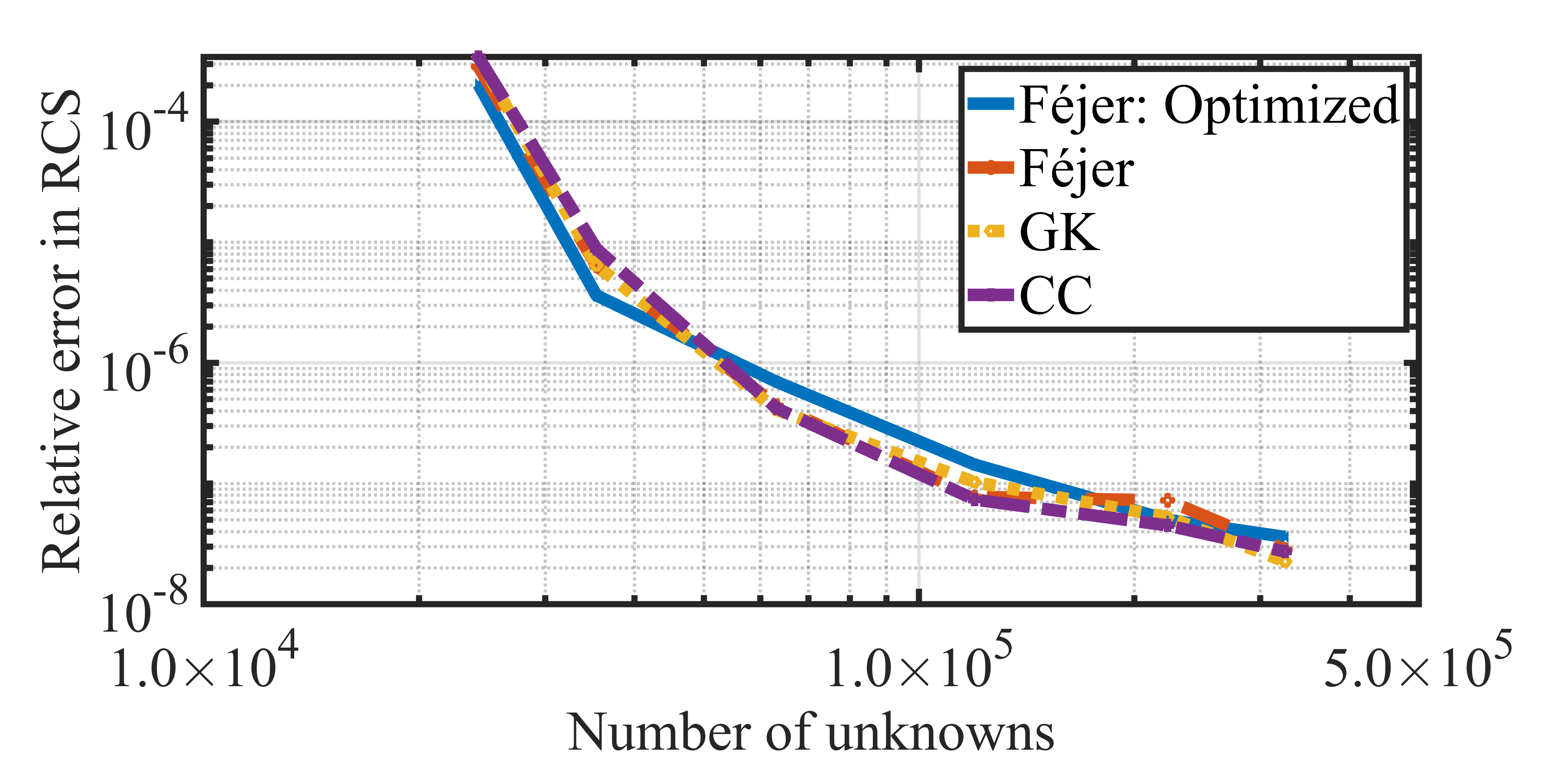}}
        
    \subfloat[]{\label{fig:sphere_rcs_vs_kernels}\includegraphics[width=1.0\columnwidth]{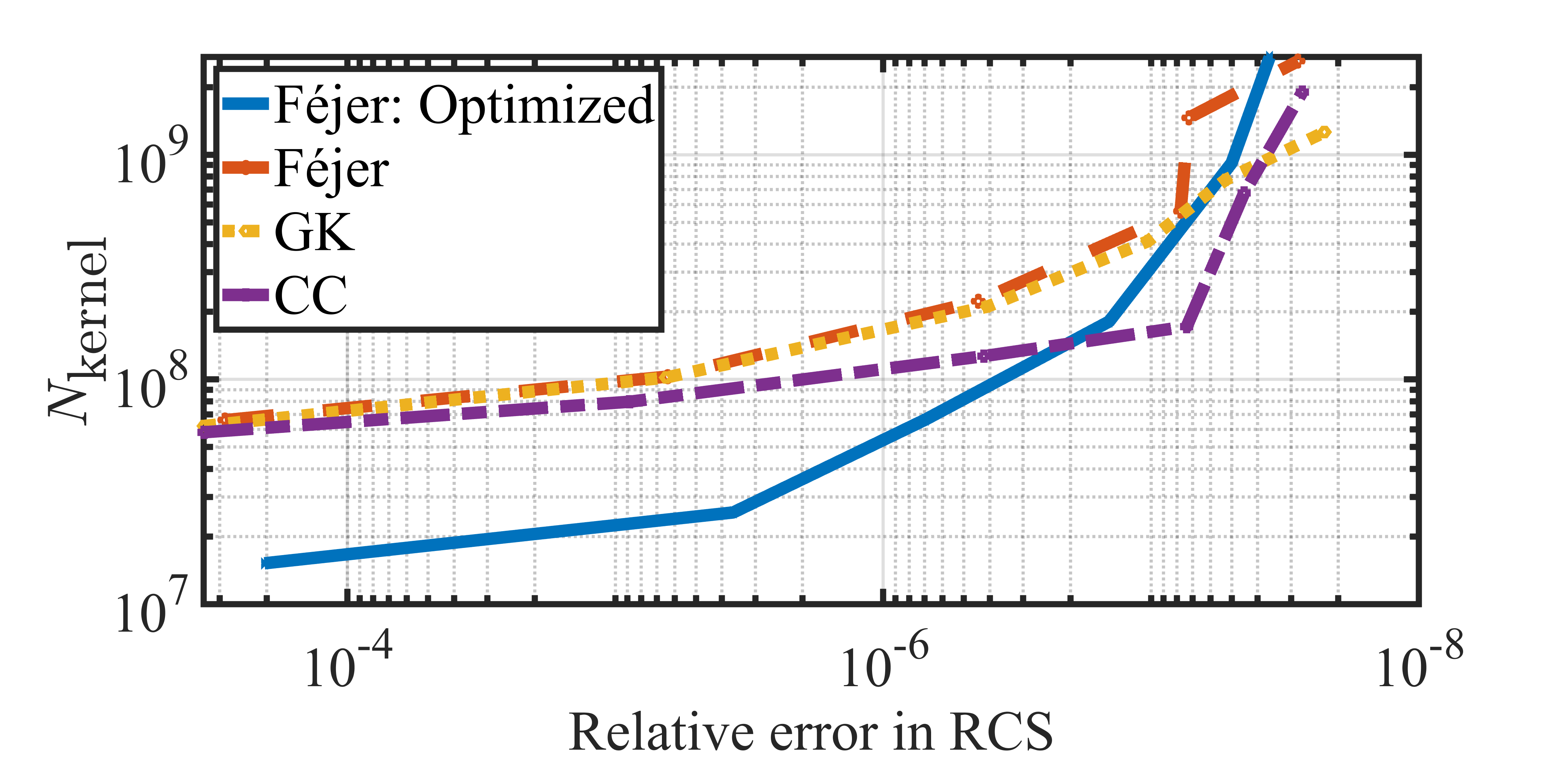}}
    
    \subfloat[]{\label{fig:sphere_rcs_vs_time}\includegraphics[width=1.0\columnwidth]{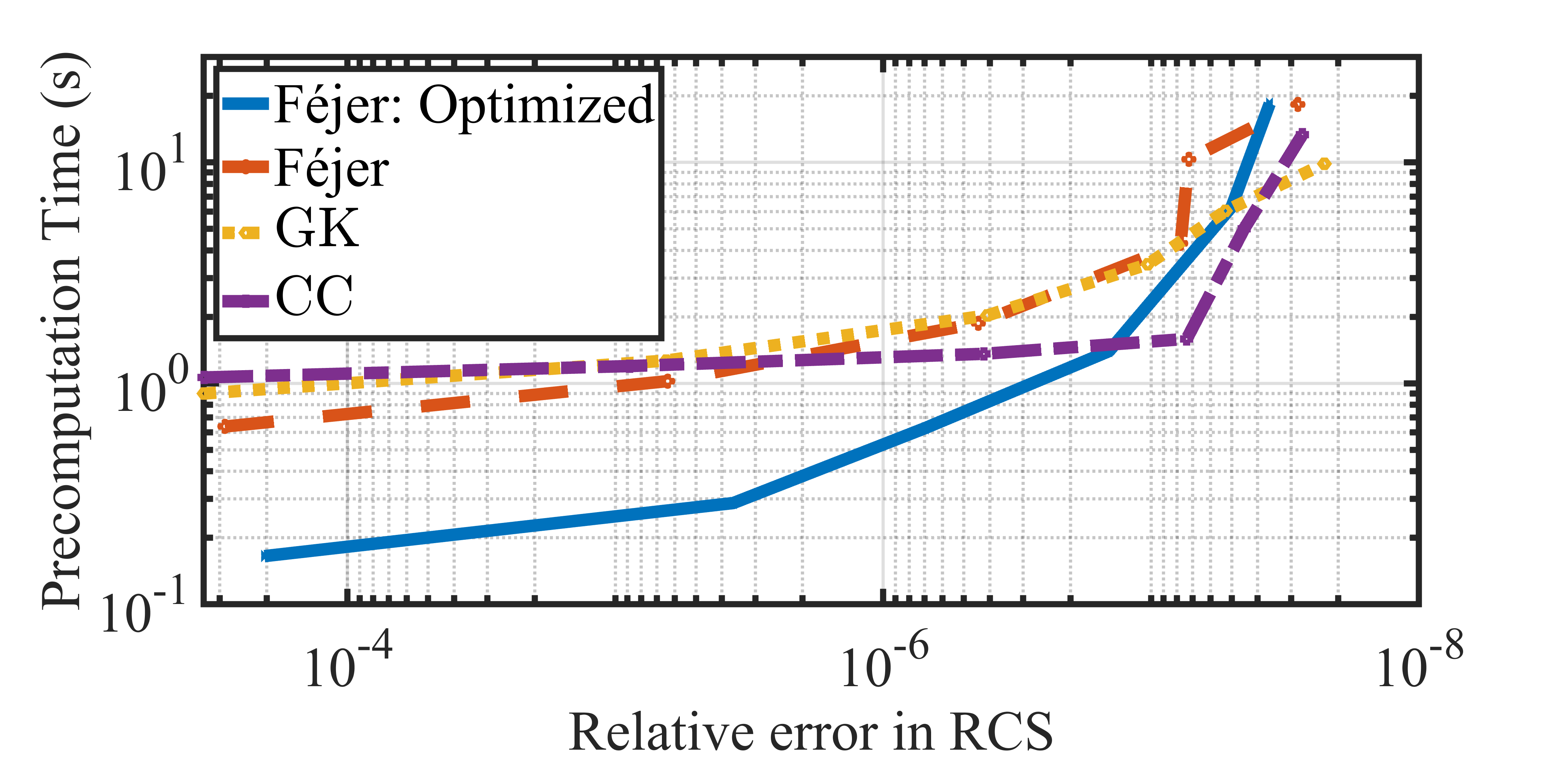}}
    
    \subfloat[]{\label{fig:sphere_rcs_vs_mem}\includegraphics[width=1.0\columnwidth]{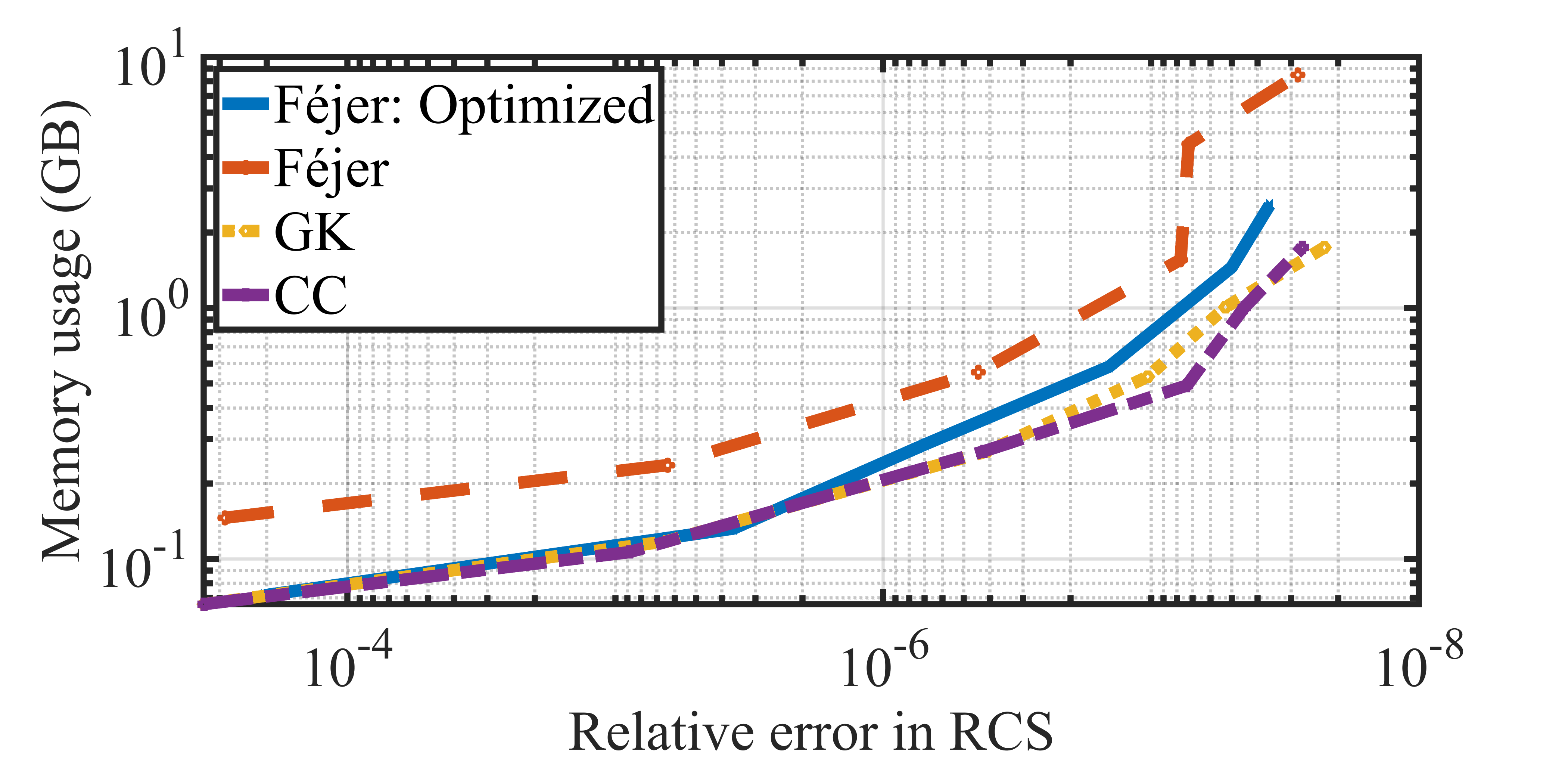}}

    \caption{Scattering by a non-uniformly split PEC Sphere of $D = 4\lambda$ diameter. Convergence of the Radar Cross Section measurement at $\theta \in [0^{\circ},180^{\circ}], \varphi = 90^{\circ}$ against (a) the number of unknowns, (b) the number of kernel evaluations, (c) the precomputation time, (d) the precomputation memory usage.}
\label{fig:sphere_rcs_parent}
\end{figure}

\subsection{Clenshaw--Curtis Quadrature}\label{CCsection}
The $p$-adaptive strategy is based on F\'ejer's second quadrature rule, and the same nodes as Clenshaw-Curtis quadrature with the exception of the endpoints~\cite{clenshaw1960method}. Quadrature nodes are arranged such that every other node of a  $2N$ order rule coincides with the nodes of an $N$ order rule. This allows for a development of an adaptive integration strategy, where the nodes of the higher order quadrature rule reuse the nodes of the lower order rule. In one dimension, the integral is approximated with an order $N$ ($I_1$) rule, and then with an order $2N$ ($I_2$) rule. The error estimate $|I_2 - I_1|$ is compared against the tolerance condition:
\[
|I_2 - I_1| < \max(\text{tol}_{\text{abs}}, \text{tol}_{\text{rel}} \cdot |I_2|).
\]
If the error tolerances are met, the result $I_2$ of higher-order rule is chosen. Otherwise, the algorithm repeatedly refines the quadrature rule (doubling the number of nodes at each step such that there are $2^{k-1}N$ points for the $k$-th iteration) and comparing against the result from the previous iteration until the tolerance is reached or the maximum refinement depth $p$ is exceeded.
The integral $I_{k}$ at the $k$-th iteration is defined as
\[
I_{k}(\mathbf{r}) \approx \sum_{i=0}^{n-1}\sum_{j=0}^{m-1} a_{i,j}\beta_{i,j}^{k},
\]
where $a_{i,j}$ are the Chebyshev expansion coefficients for the ``auxiliary" density, and $\beta_{i,j}^{k}$ denotes the result of the precomputation integral ~\eqref{eq:nearinteraction_weights}.

\begin{figure}[!ht]
    \centering
    \subfloat[]{\label{fig:sphere_mie_vs_unknowns}\includegraphics[width=1.0\columnwidth]{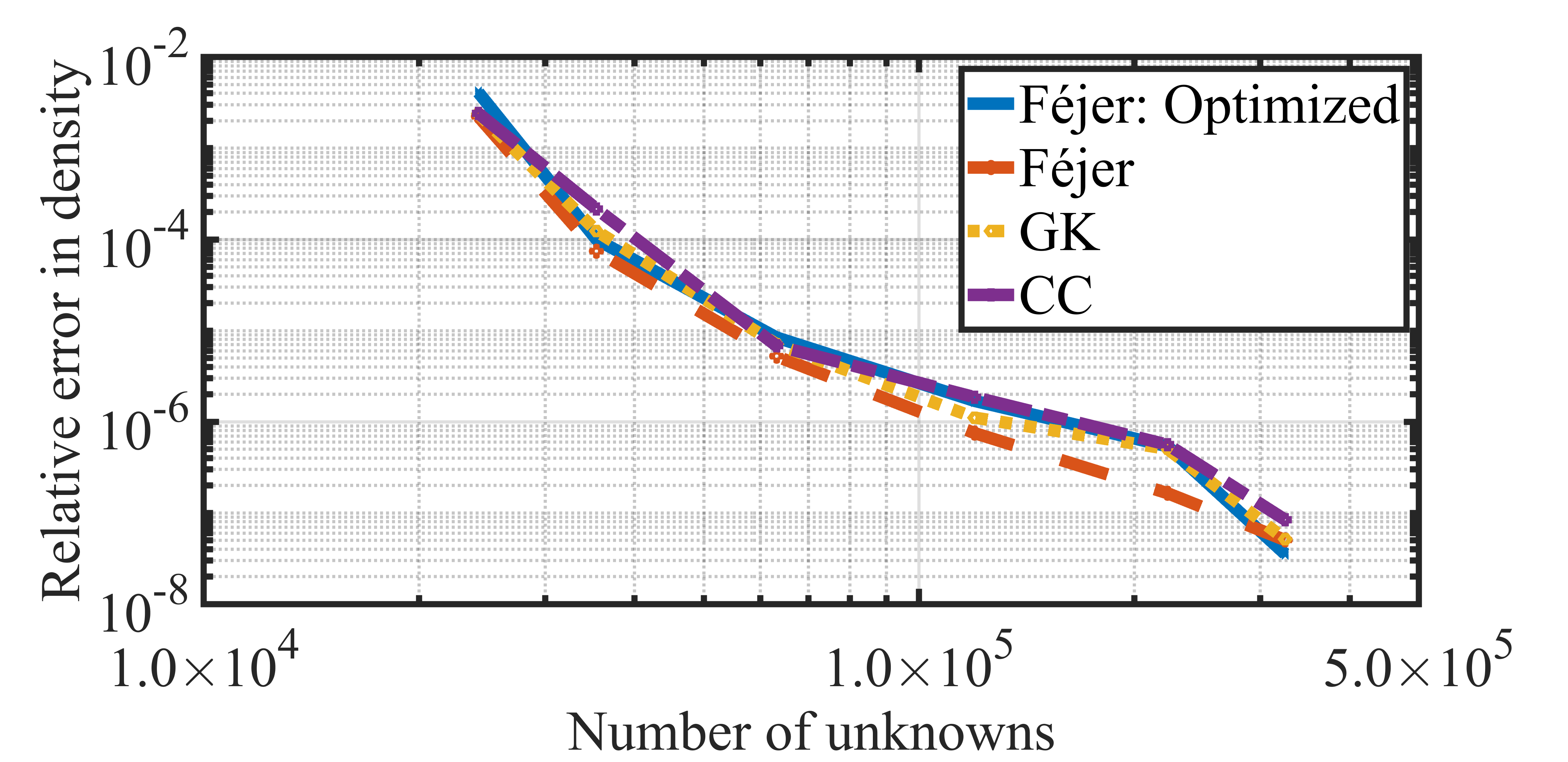}}
    
    \subfloat[]{\label{fig:sphere_mie_vs_kernels}\includegraphics[width=1.0\columnwidth]{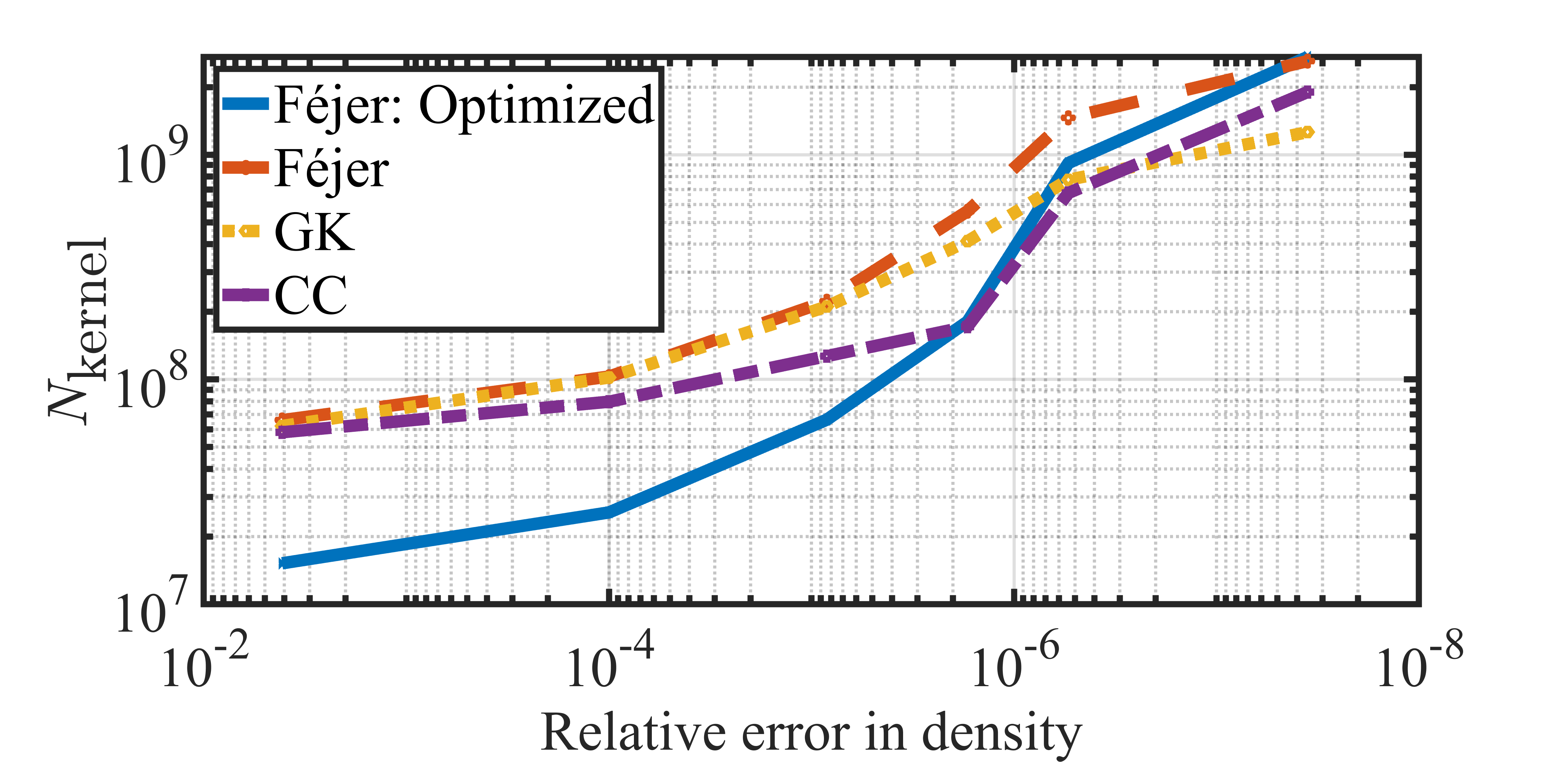}}
    
    \subfloat[]{\label{fig:sphere_mie_vs_time}\includegraphics[width=1.0\columnwidth]{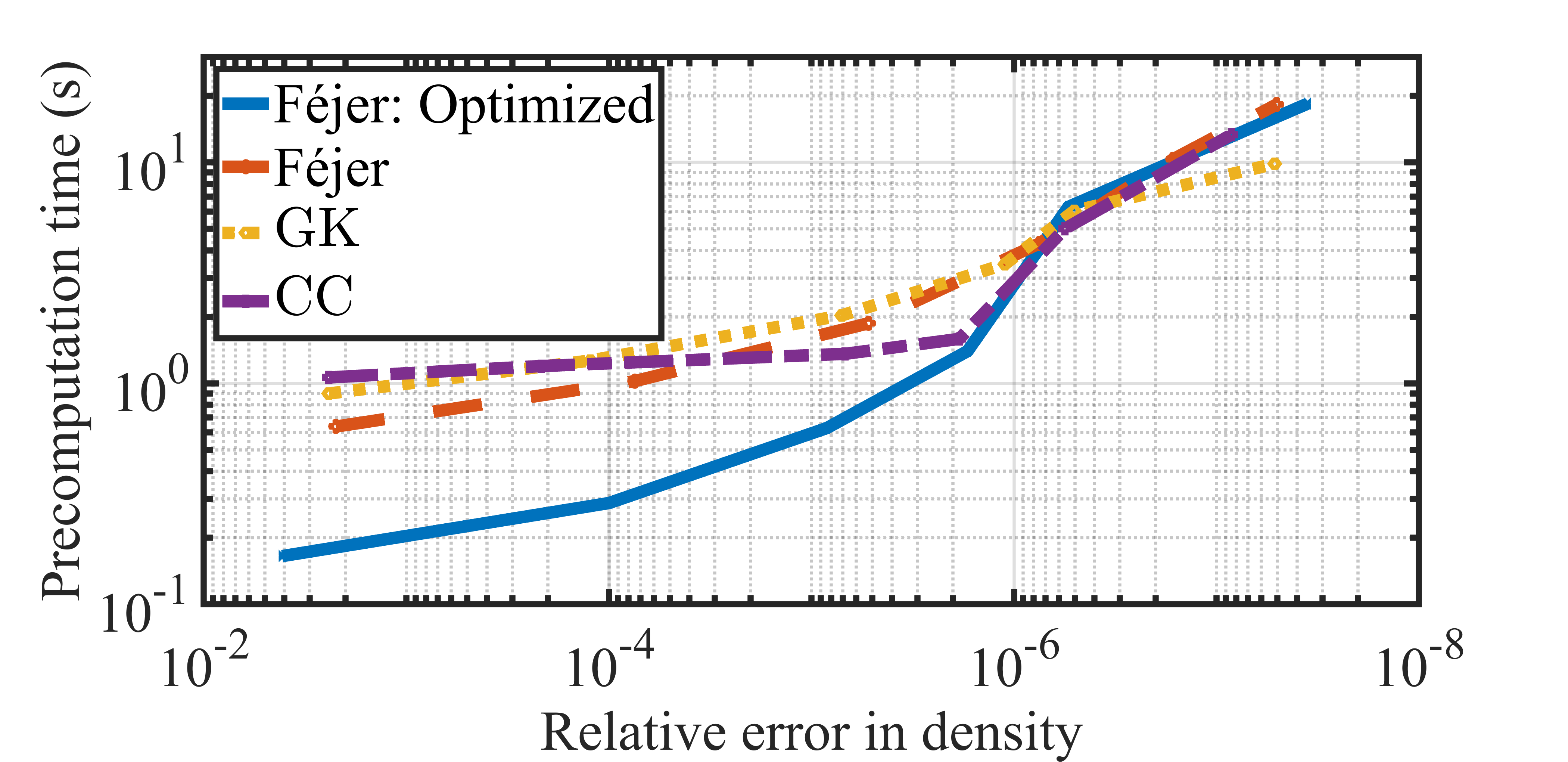}}
    
    \subfloat[]{\label{fig:sphere_mie_vs_mem}\includegraphics[width=1.0\columnwidth]{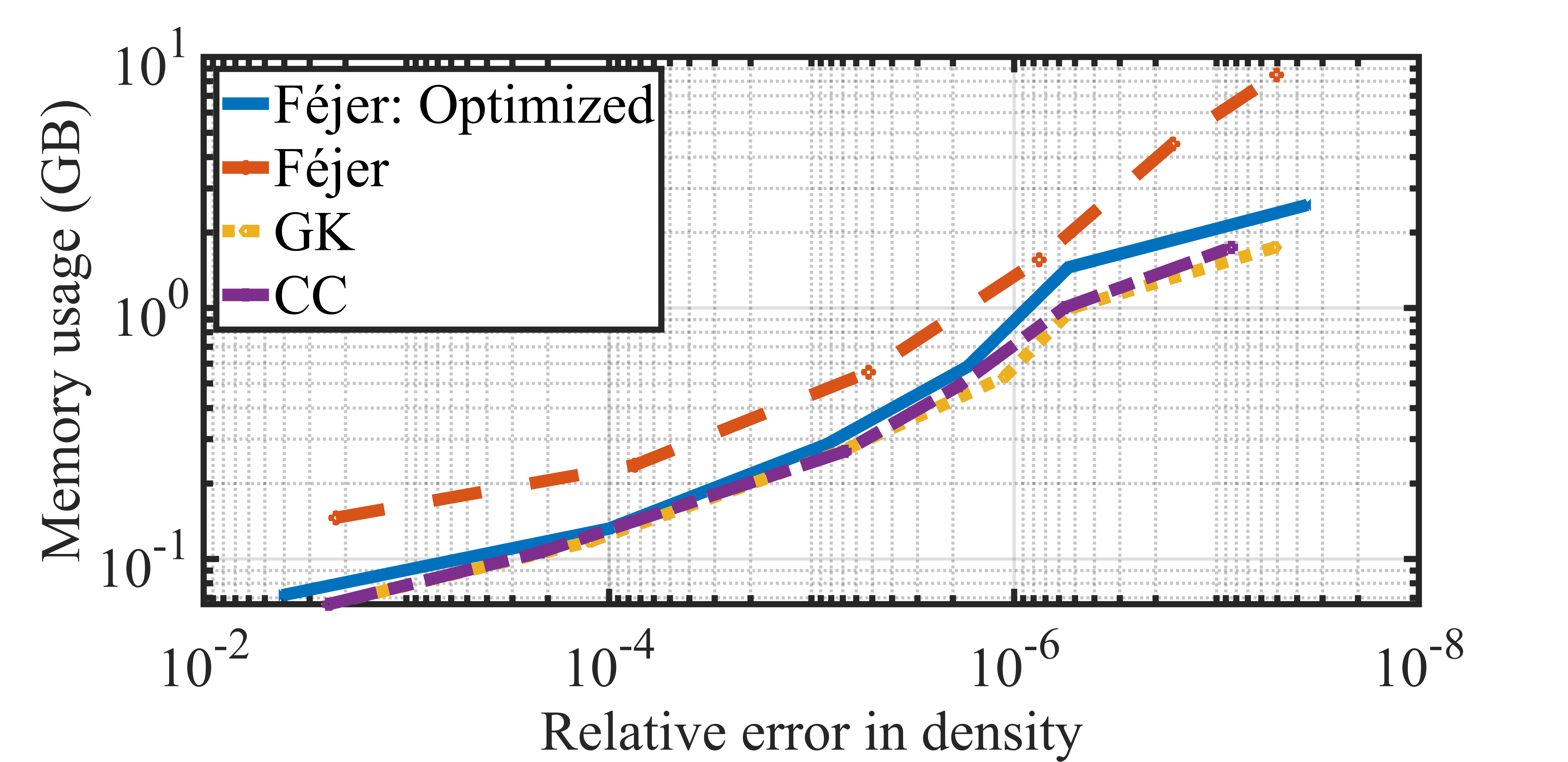}}

    \caption{Scattering by a non-uniformly split PEC Sphere of $D = 4\lambda$ diameter. Convergence of the surface current density against (a) the number of unknowns, (b) the number of kernel evaluations, (c) the precomputation time, (d) the precomputation memory usage.}
\label{fig:mie_parent}
\end{figure}

\begin{figure}[!t]
\centerline{\includegraphics[width=0.75\columnwidth]{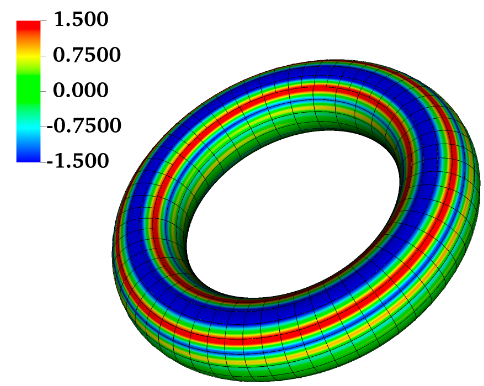}}
\caption{The real part of $J_x$ solved surface current density on a PEC toroid of $D = 26.2\lambda$ ``outside" and $d = 15.7\lambda$ ``inside" diameters.}
\label{fig:toroid_density}
\end{figure}

\begin{figure}[!t]
    \centering

    \subfloat[]{\label{fig:toroid_rcs_vs_unknowns}\includegraphics[width=1.0\columnwidth]{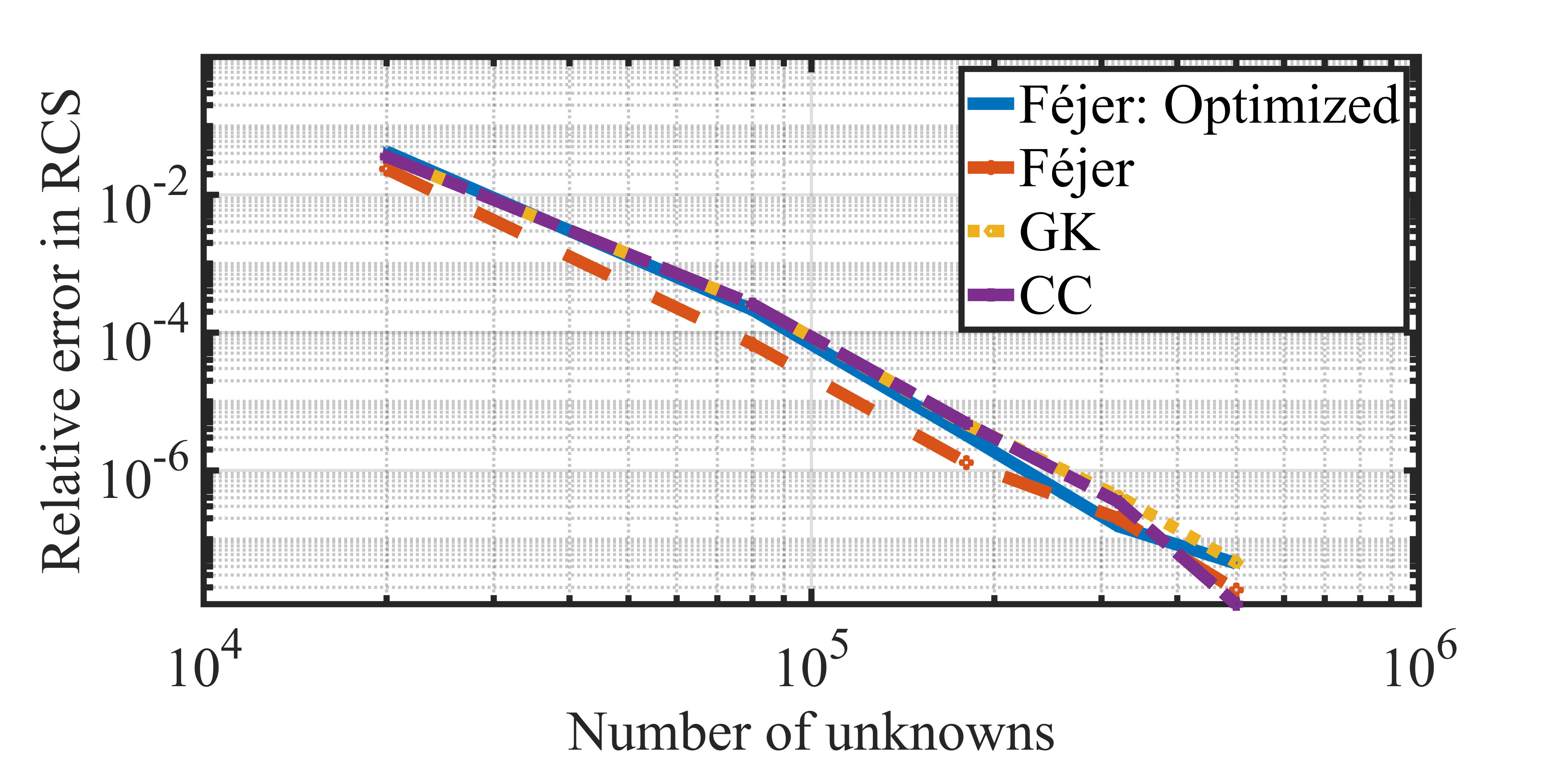}}
    
    \subfloat[]{\label{fig:toroid_rcs_vs_kernels}\includegraphics[width=1.0\columnwidth]{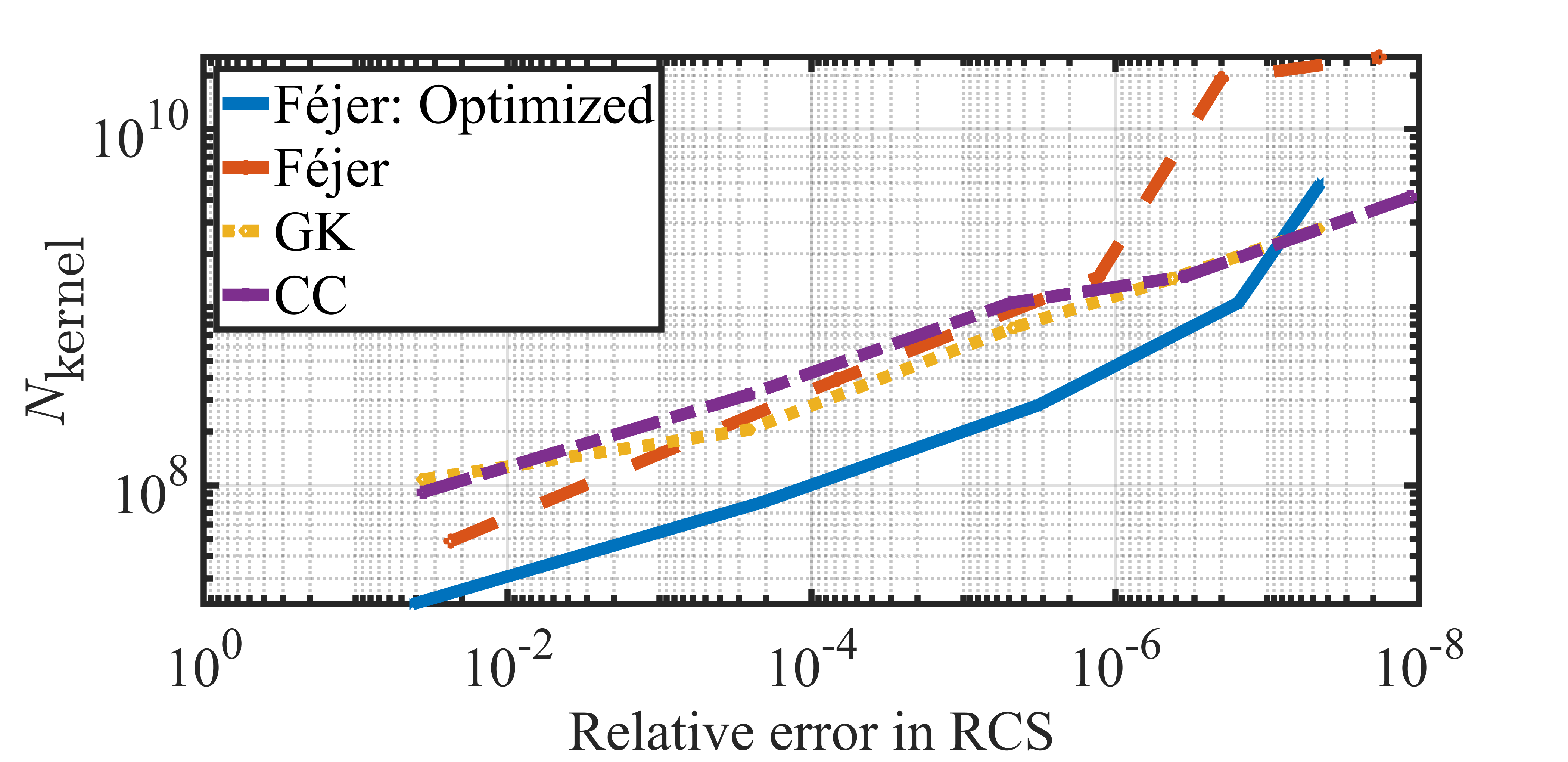}}
    
    \subfloat[]{\label{fig:toroid_rcs_vs_time}\includegraphics[width=1.0\columnwidth]{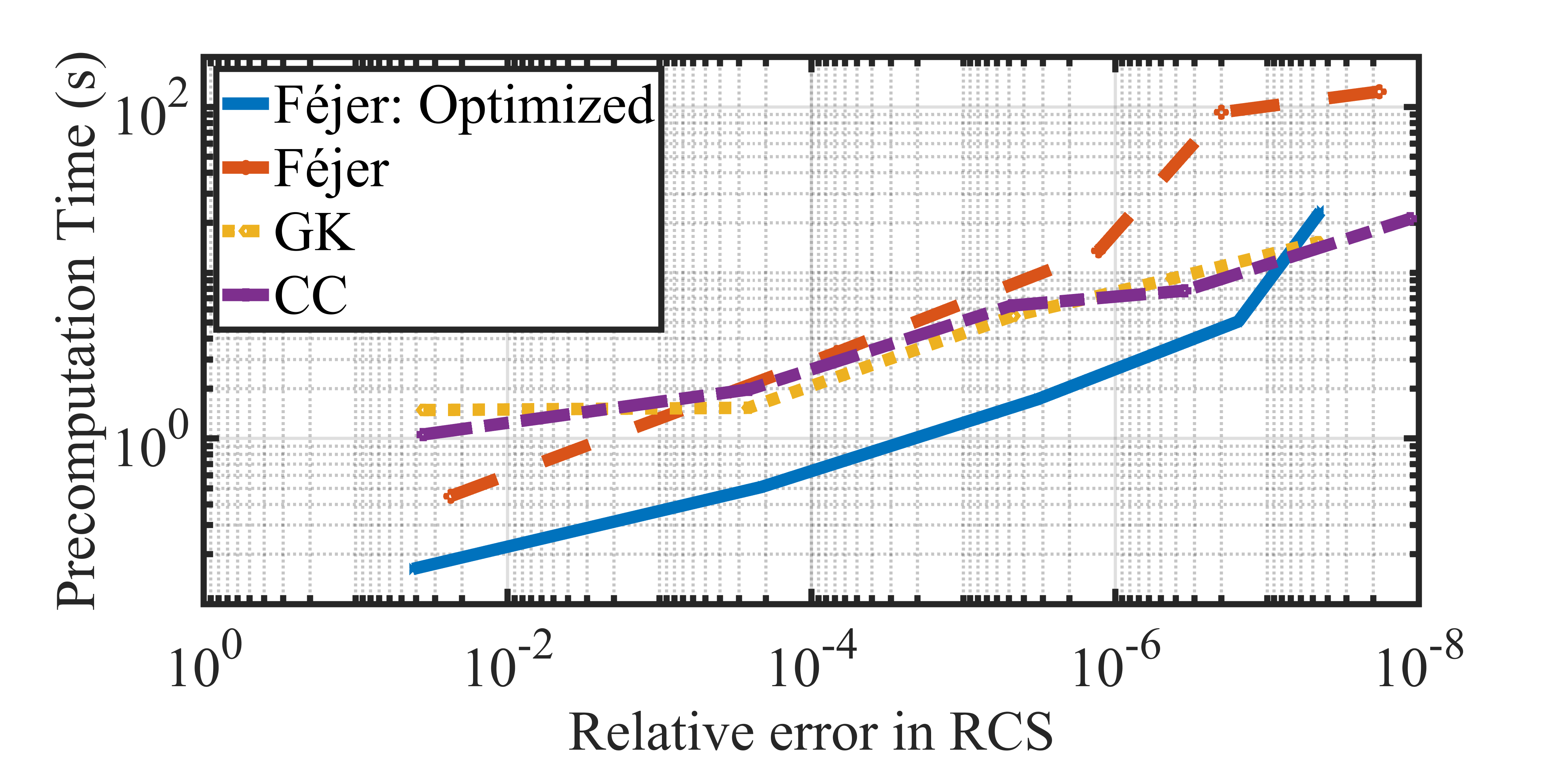}}
    
    \subfloat[]{\label{fig:toroid_rcs_vs_mem}\includegraphics[width=1.0\columnwidth]{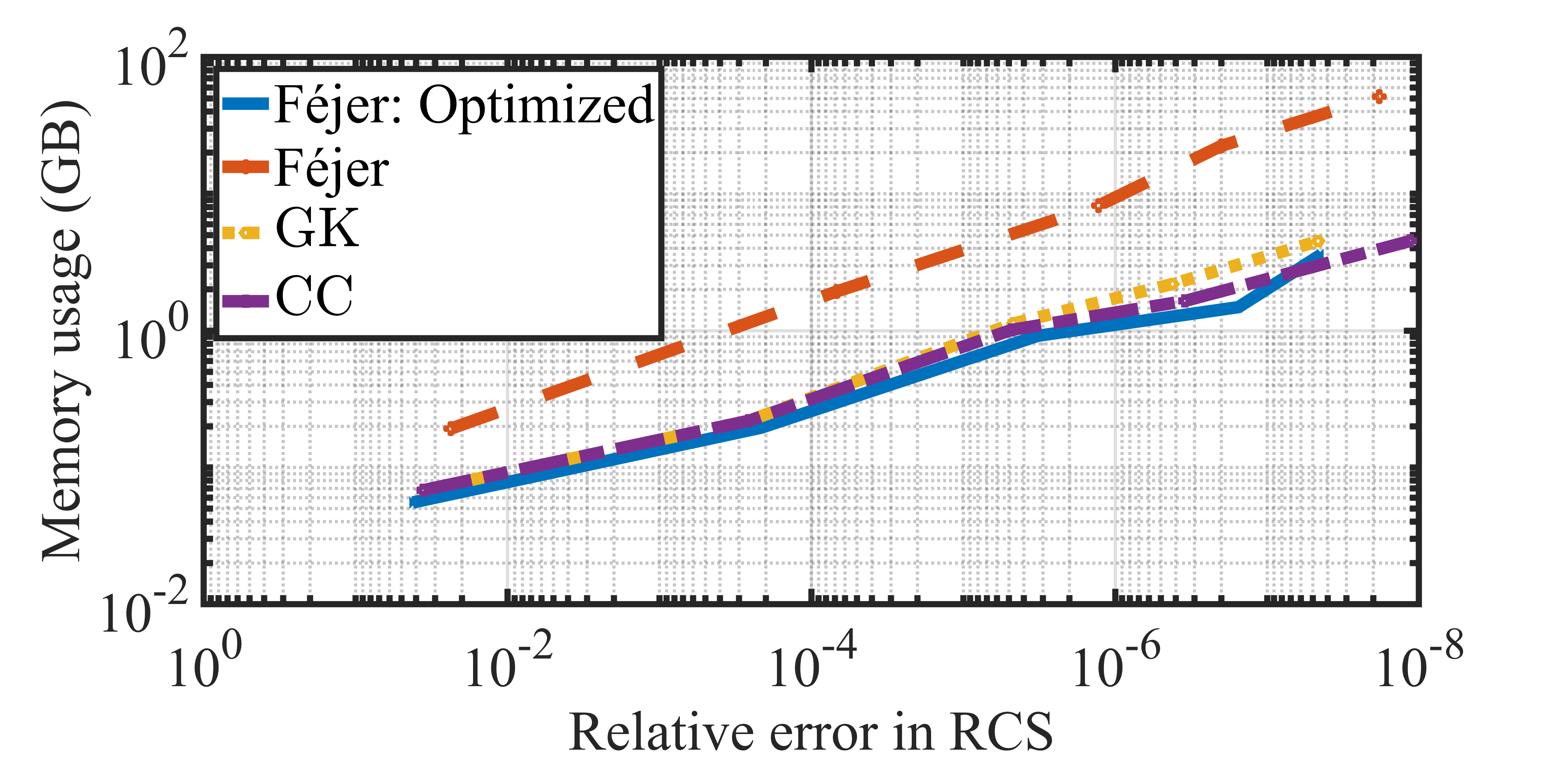}}

    \caption{Scattering by a uniformly split PEC Toroid of $D = 26.2\lambda$ ``outside" diameter and $d = 15.7\lambda$ ``inner" diameter. Convergence of the Radar Cross Section measurement at $\theta \in [0^{\circ},180^{\circ}], \varphi = 90^{\circ}$ against (a) the number of unknowns, (b) the number of kernel evaluations, (c) the precomputation time, (d) the precomputation memory usage.}
\label{fig:toroid_parent}
\end{figure}

\begin{figure}[!t]
\centerline{\includegraphics[width=0.8\columnwidth]{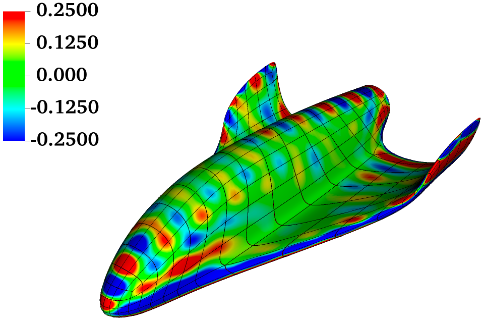}}
\caption{The real part of $J_x$ solved surface current density on PEC Glider of $W = 7.5\lambda$ wingspan and $L = 10\lambda$ length.}
\label{fig:glider_density}
\end{figure}

\begin{figure}[!ht]
    \centering

    \subfloat[]{\label{fig:glider_rcs_vs_unknowns}\includegraphics[width=1.0\columnwidth]{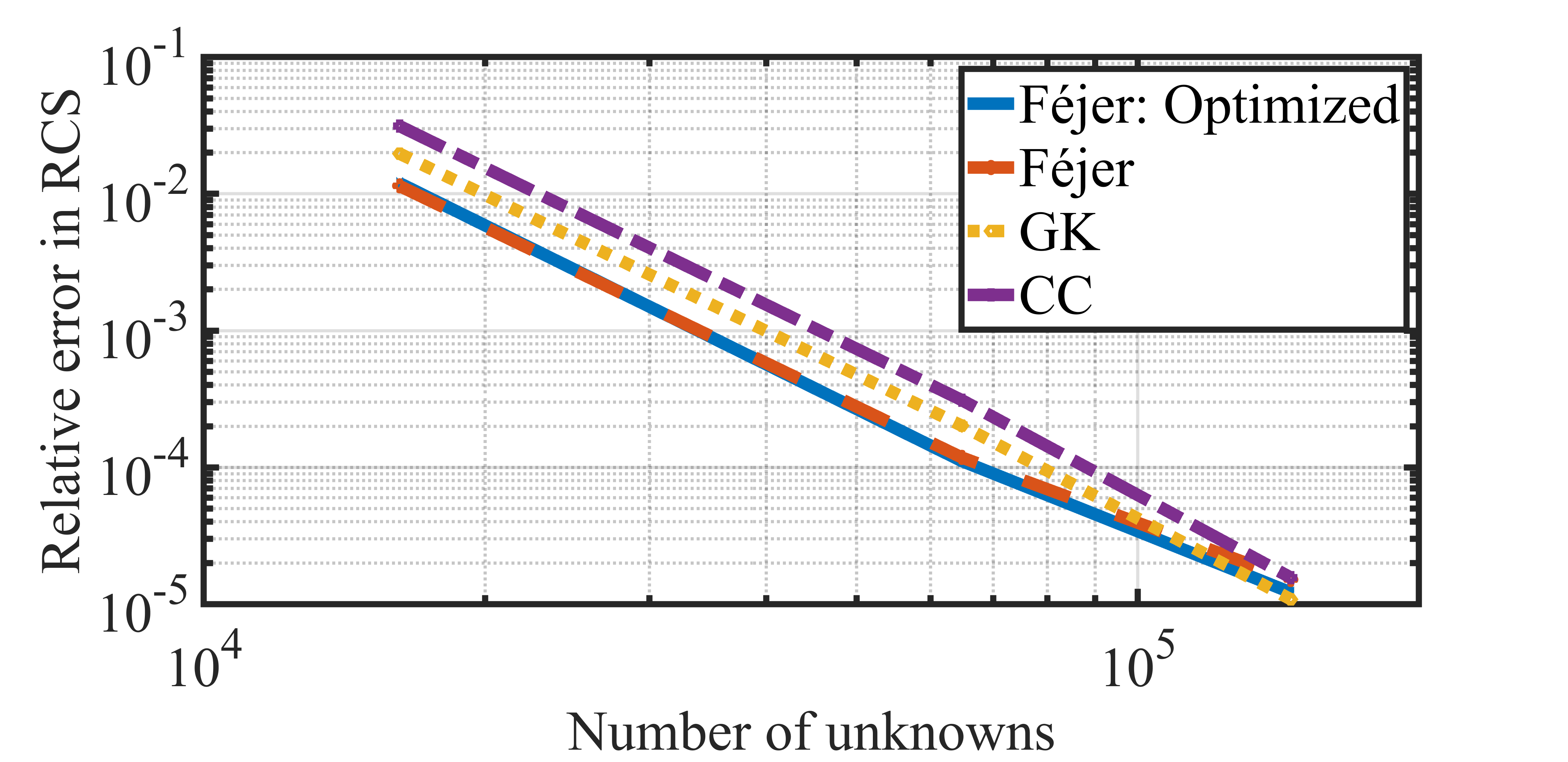}}
    
    \subfloat[]{\label{fig:glider_rcs_vs_kernels}\includegraphics[width=1.0\columnwidth]{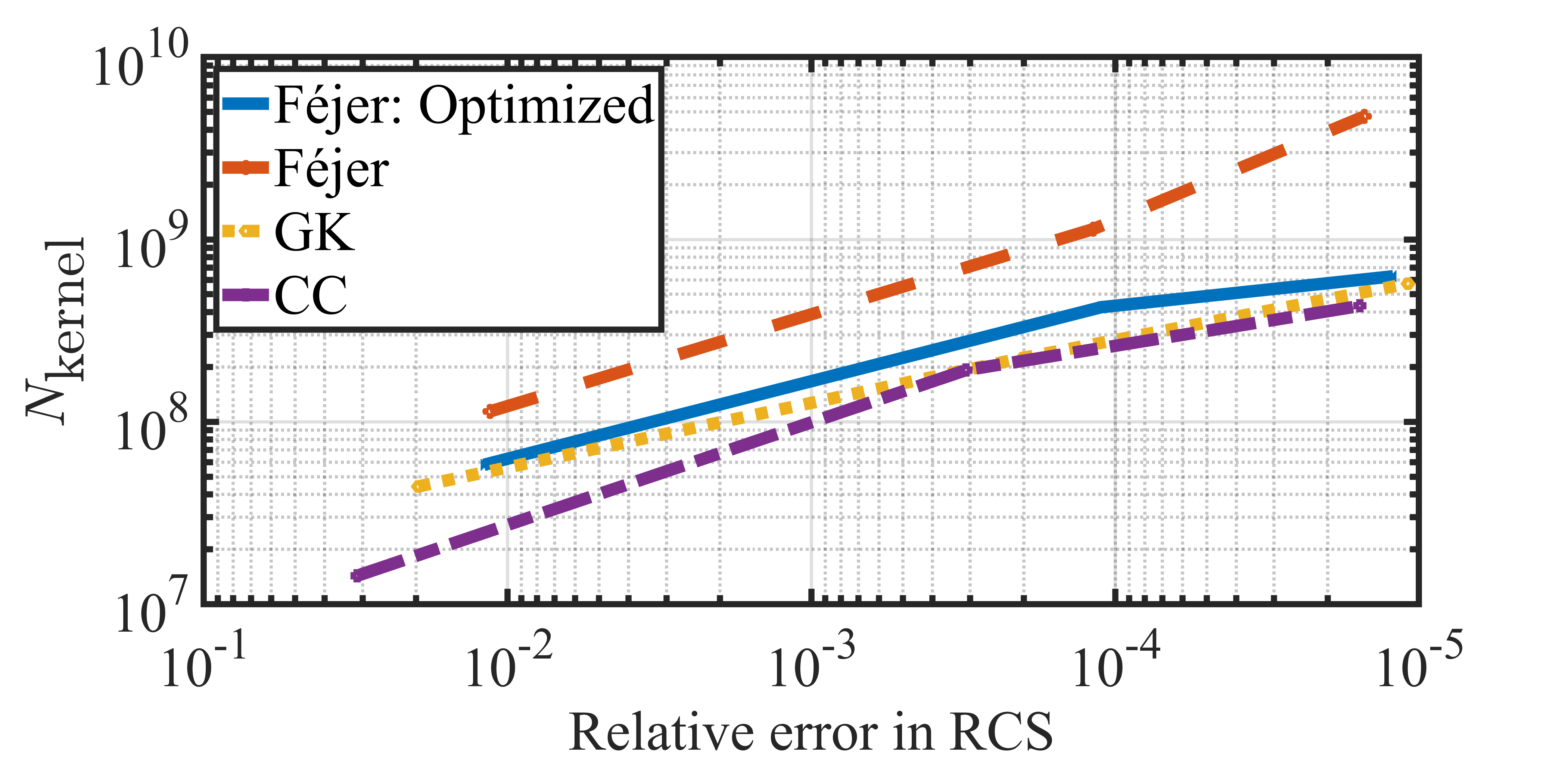}}
    
    \subfloat[]{\label{fig:glider_rcs_vs_time}\includegraphics[width=1.0\columnwidth]{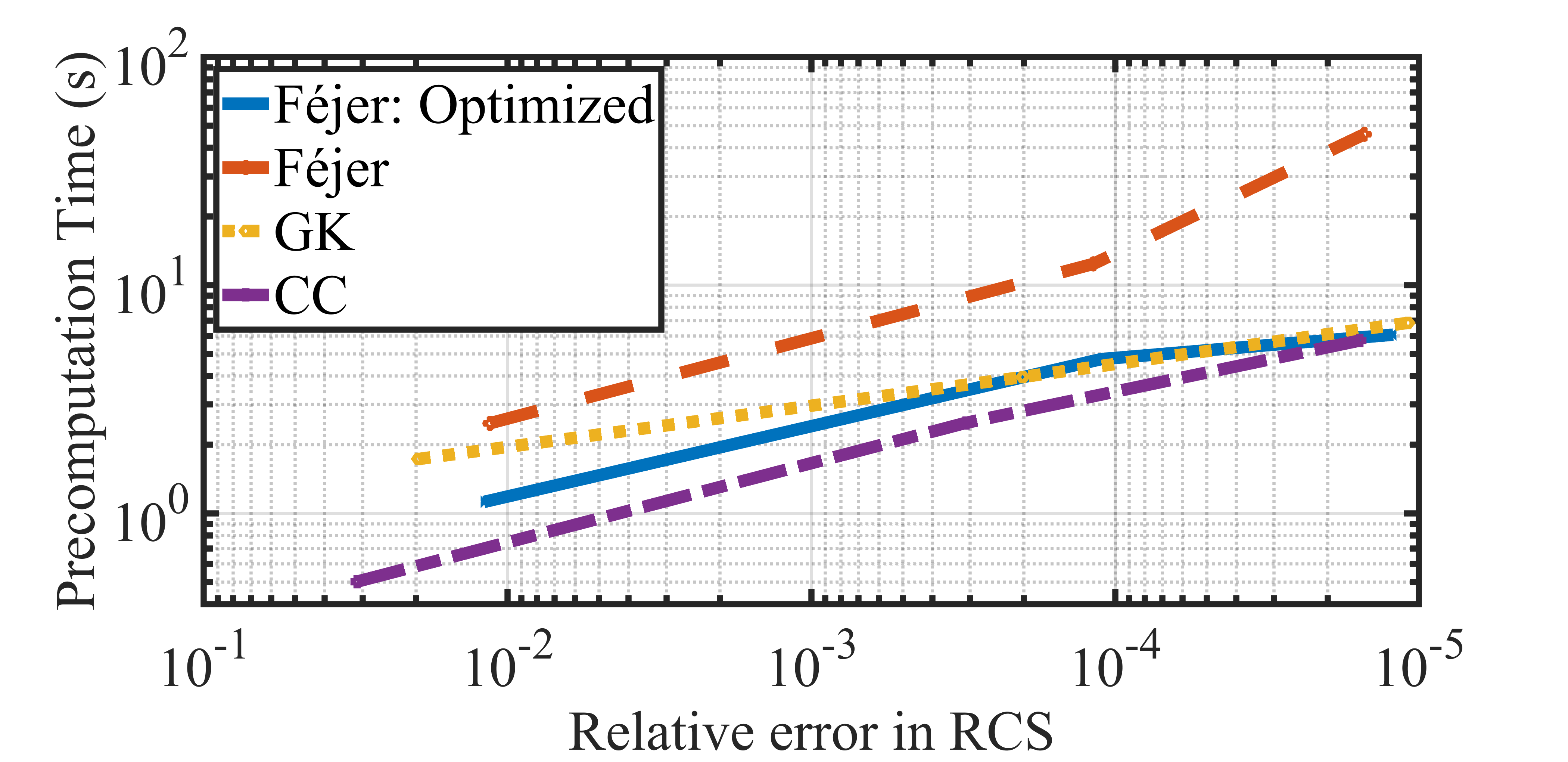}}
    
    \subfloat[]{\label{fig:glider_rcs_vs_mem}\includegraphics[width=1.0\columnwidth]{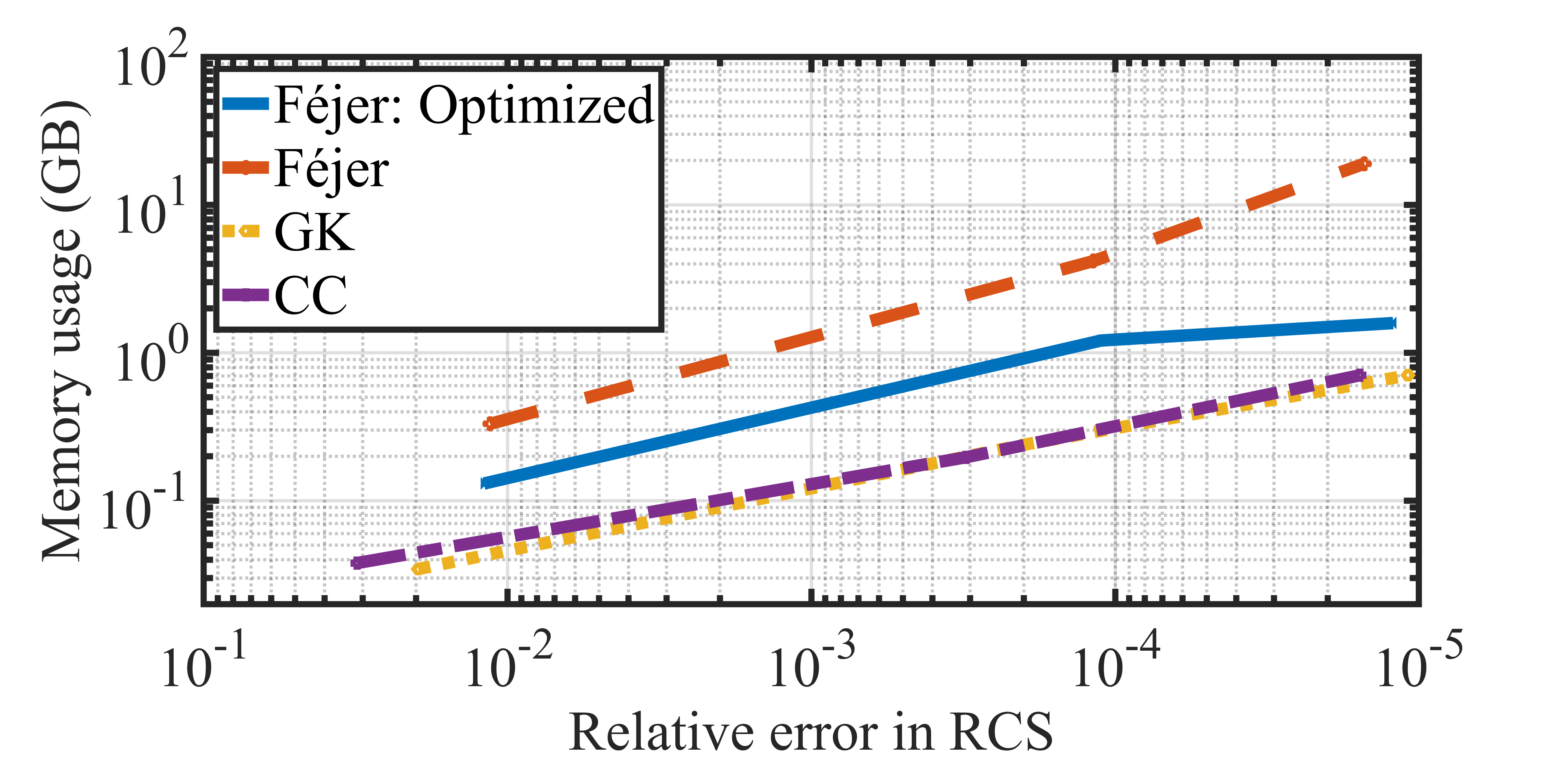}}

    \caption{Scattering by a PEC CAD glider model of $L = 10\lambda$ wingspan and $W = 7.5\lambda$ length. Convergence of the Radar Cross Section measurement at $\theta \in [0^{\circ},180^{\circ}], \varphi = 90^{\circ}$ against (a) the number of unknowns, (b) the number of kernel evaluations, (c) the precomputation time, (d) the precomputation memory usage.}
\label{fig:glider_parent}
\end{figure}

\section{Numerical Results}\label{results}

Many real-world electromagnetic scattering problems involve geometrically complex objects that cannot be discretized using uniformly sized patches $\Gamma_p$ and might require the use of ``defective" meshes, where adjacent patches may not share edges perfectly. The Nyström method can easily handle such geometric representations. This section presents numerical validation of the proposed adaptive integration approach on three geometries of increasing complexity: a non-uniformly discretized PEC sphere, a uniformly meshed toroidal scatterer, and a CAD model of a glider. An incident plane wave excitation with a wavenumber $k_{\mathrm{aux}} = 1.1\cdot k_0$ was used as the ``auxiliary" density for error estimation. The coefficient can be modified based on the complexity of the problem and the geometry; however, the chosen coefficient empirically showed robust performance for all three objects considered. All simulations were run on a server with dual AMD EPYC 7763 processors with 128 cores.

The first example is a PEC sphere (Fig. \ref{fig:sphere_density}) with a $D=4\lambda$ diameter and serves as a validation case. The scattering problem for a spherical scatterer has a closed-form analytical solution in the form of the Mie series, which  will be used as the reference. Figs.~\ref{fig:sphere_rcs_parent} and~\ref{fig:mie_parent} show the convergence of both the Radar Cross Section (RCS) and the surface current density with respect to the number of unknowns, kernel evaluations, the precomputation time, and the amount of memory (RAM) used to store the precomputations.

We compare the results from four strategies of calculating the precomputations: ``F\'ejer: Optimized", ``F\'ejer", ``GK", and ``CC". The ``F\'ejer: Optimized" represents a fixed-grid ``F\'ejer" quadrature setup in which both the singular refinement grid and near-interaction threshold $\Delta_{\mathrm{near}}$ have been meticulously hand tuned across multiple simulation runs. This represents the most time-efficient implementation of the fixed quadrature approach, as it employs the coarsest grid and smallest $\Delta_{\mathrm{near}}$ that maintain accuracy within 0.5 digits of the best accuracy possible for a given discretization. The standard ``F\'ejer" uses the same fixed quadrature settings and $\Delta_{\mathrm{near}}$ across all discretizations. These settings were optimized for a single discretization (where the surface is represented by 1624 quadrilateral patches and 324800 unknowns and is the most refined discretization considered for the scatterer) and then reused, representing a more typical usage scenario in practice. While the standard ``F\'ejer" approach is expected to match the accuracy of the ``Optimized" version, it is expected to exhibit lower computational efficiency. {GK} and {CC} denote results obtained using the proposed $h$-adaptive Gauss-Kronrod and $p$-adaptive Clenshaw-Curtis integration methods, respectively, as described in Sections~\ref{GKsection} and~\ref{CCsection}.

The results show that adaptive approaches ``GK" and ``CC" achieve better than $4 \cdot 10^{-8}$ relative error in RCS and surface density, and perform no worse than the ``F\'ejer: Optimized'' scenario. This is especially significant since the latter assumes prior knowledge of optimal tuning parameters, which required solving the problem many times while adjusting the parameters to optimize the performance. This approach is completely impractical for realistic scenarios. Compared to the realistic use case described by the ``F\'ejer'' baseline, the adaptive methods demonstrate savings in the number of singular kernel evaluations, precomputation time, and memory usage, without compromising any accuracy.

Next, we consider a toroid scatterer (\ref{fig:toroid_density}) with an ``external" diameter of $D=26.2\lambda$ and an ``internal" diameter of $d=15.7\lambda$. Fig.~\ref {fig:toroid_parent} illustrates RCS convergence for the toroidal PEC scatterer. This geometry lacks an analytical solution, so a highly discretized numerical solution is used as a reference. For the toroid scatter, ``GK" and ``CC" both exhibit strong performance. They achieve a similar accuracy as the fixed-grid methods with comparable (for ``F\'ejer: Optimized") or better (for the more practical use-case of ``F\'ejer") wall-clock time. Although the meticulously optimized fixed-grid method occasionally achieves a comparable accuracy with a marginally better wall clock time (on the order of 1-2 seconds), the adaptive methods provide better generalization, achieve the desired accuracy with a single run, and do not require problem-specific parameter tuning.

\begin{figure}[!h]
\centerline{\includegraphics[width=1.0\columnwidth]{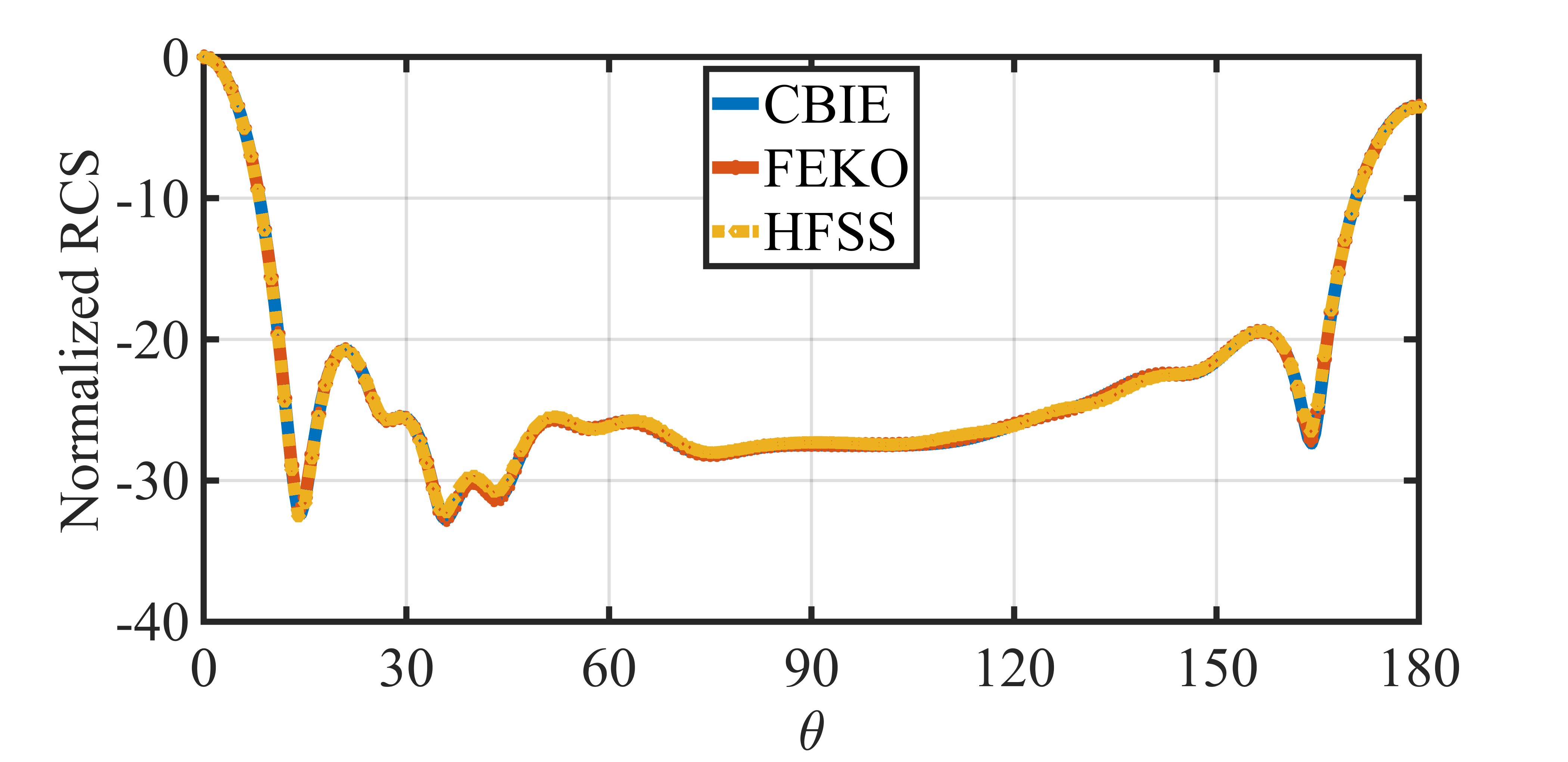}}
\caption{RCS at $\phi = 90^{\circ}$ for the PEC glider model calculated using Ansys HFSS, Altair FEKO and the Adaptive CBIE method.}
\label{fig:hfss_feko_glider}
\end{figure}

The third example involves a CAD model of a glider with a $W=7.5\lambda$ wingspan and $L=10\lambda$ length. The model has pronounced geometric features, including sharp edges near the exhaust. As shown in Fig.~\ref{fig:glider_parent}, the adaptive methods maintain robust convergence, and similar accuracy to the fixed-grid methods. The manually parameter tuned fixed-grid method occasionally achieves a comparable accuracy in a similar wall clock time; however, for complex cases the adaptive methods can provide up to 8X savings in wall clock time for realistic use-cases. The accuracy of all methods is limited to approximately 5 digits, which is attributed to the sharp-edge singularities that are challenging for high-order methods. Future work will address these limitations by regularizing the integral equation operators and leveraging the corner-singularity handling techniques presented in \cite{2dbie, aces_2dbie}. We validate our results by comparing the RCS pattern against the $\theta$ at $\phi = 90^{\circ}$ calculated using the CBIE solver with adaptive integration approaches against the commercial software Ansys HFSS and Altair FEKO in Fig. \ref{fig:hfss_feko_glider}. The patterns produced by all three solvers match very closely.

\begin{figure}[h]
    \centering

    \subfloat[]{\label{fig:sphere_near_dist}\includegraphics[width=0.75\columnwidth]{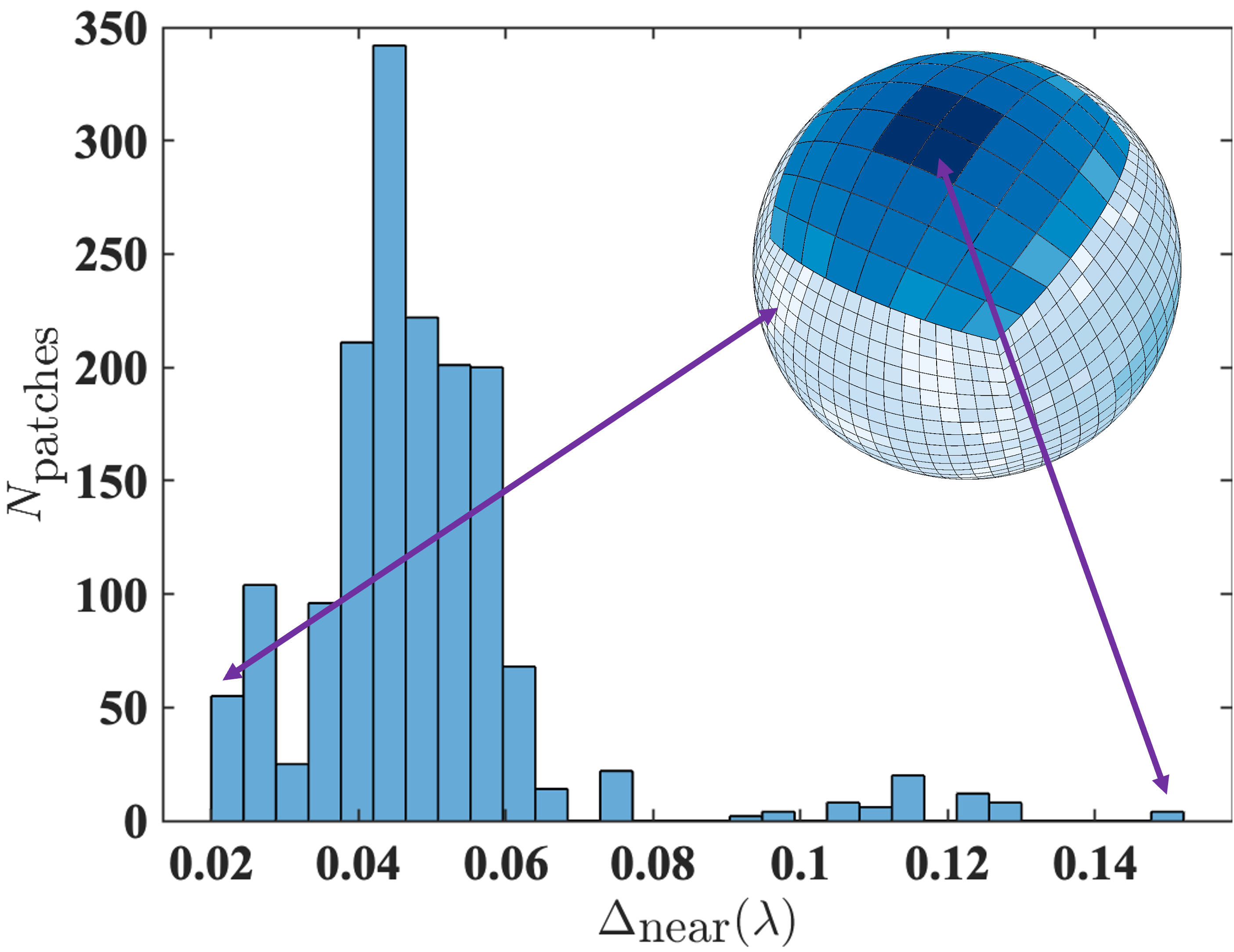}}
    
    \subfloat[]{\label{fig:toroid_near_dist}\includegraphics[width=0.75\columnwidth]{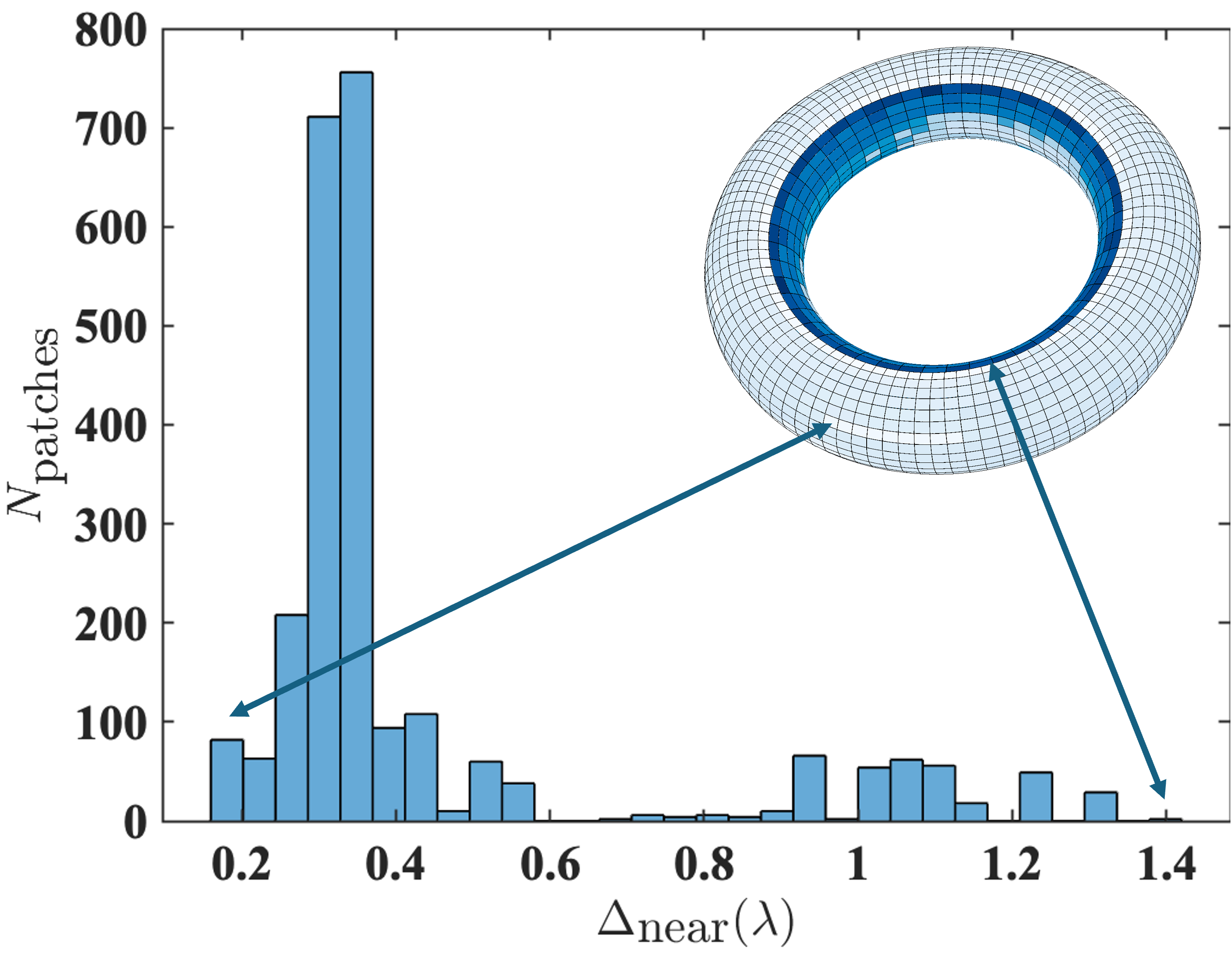}}
    
    \subfloat[]{\label{fig:glider_near_dist}\includegraphics[width=0.75\columnwidth]{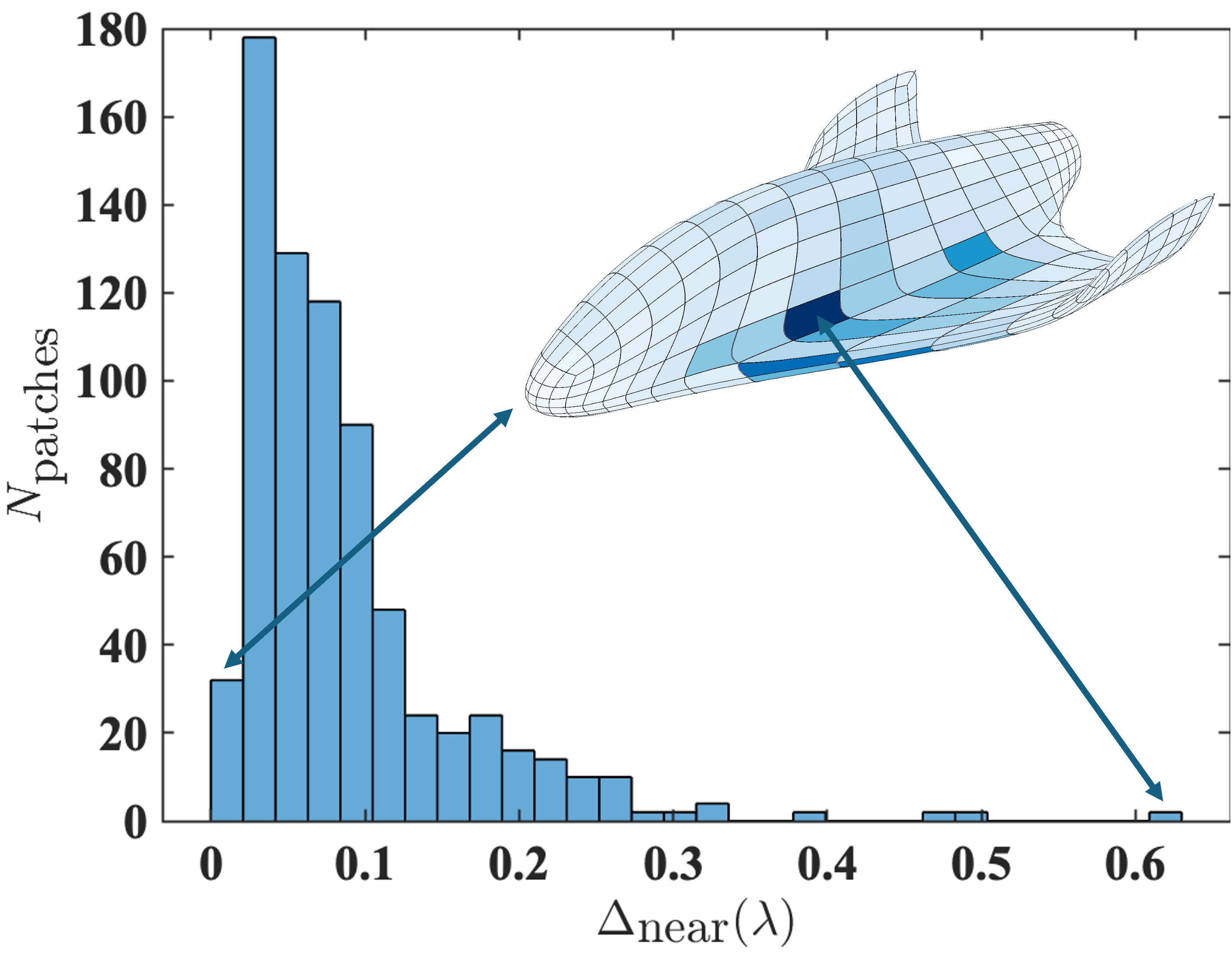}}
    
    \caption{$\Delta_{\textrm{near}}$ parameter distribution for integration patches for (a) a PEC Sphere, (b) a PEC Toroid, (c) a PEC Glider scatterers.}
\label{fig:near_dist_parent}
\end{figure}

Fig. \ref{fig:near_dist_parent} describes the number of integration patches per $\Delta_{\textrm{near}}$ parameter across the considered geometries. $\Delta_{\textrm{near}}$ is chosen by the adaptive algorithms and is in units of the wavelength $\lambda$. The darker colors in the colormap represent larger $\Delta_{\textrm{near}}$ values. For these patches, more target points are considered ``near-singular" and require the use of singular integration approach to achieve the desired accuracy. This requires higher memory usage for storing the precomputation weights. 

Fig. \ref{fig:sphere_near_dist} shows the distribution for the sphere parametrization with 1624 patches. The majority of integration patches require small $\Delta_{\textrm{near}}$ values, clustered tightly around the value of $0.04\lambda$. This can be explained by the $C^{\infty}$ continuity of the object, the constant Gaussian curvature of the object, and the large number of patches used to represent the geometry. The few outlying values represent patches that are significantly larger than the rest (due to the use of defective meshes). The fixed-grid F\'ejer implementations use the same $\Delta_{\textrm{near}}$ parameter for all integration patches, which, if chosen wrong, can limit the accuracy (if the chosen $\Delta_{\textrm{near}}$ is too small) or be inefficient (if the chosen  $\Delta_{\textrm{near}}$ is too large). The adaptive methods automatically select the necessary $\Delta_{\textrm{near}}$ for each patch and can save up to 4X memory usage for the PEC Sphere scatterer (Figs. \ref{fig:sphere_rcs_vs_mem} and \ref{fig:sphere_mie_vs_mem}) when compared against the realistic use case of fixed F\'ejer method.

Fig. \ref{fig:toroid_near_dist} shows the distribution for the toroid parametrization with 2500 patches. For this parametrization, $\Delta_{\textrm{near}}$ takes a larger range of values, which can be explained by the changing curvature of the object. The curvature on the ``outside" perimeter of the toroid is smaller than that on the ``inside" perimeter, which explains the need for larger $\Delta_{\textrm{near}}$ values for patches on the ``inner" surface of the toroid. For this case, up to 11X memory usage savings is achieved (Fig. \ref{fig:toroid_rcs_vs_mem}) compared to the previous fixed F\'ejer method. Fig. \ref{fig:glider_near_dist} shows the distribution for the glider parametrization with 729 patches. For this complex geometry, the adaptive methods show up to 2X savings in memory usage compared to the manually parameter tuned fixed ``F\'ejer best" method, and up to 27X savings in memory usage (Fig. \ref{fig:glider_rcs_vs_mem}) compared to the more practical fixed F\'ejer method.

\section{Conclusions}\label{conclusions}
This manuscript introduces an adaptive approach that obviates the need for any manual parameter tuning in high-order Nystr\"om-based methods. We demonstrate the techniques by incorporating them into the Chebyshev-based discretization method and solving Maxwell's equations for 3D metallic scatterers. The approach uses  Gauss–Kronrod or Clenshaw–Curtis adaptive quadrature rules with singularity-canceling changes of variables to automatically achieve the desired accuracy in the solution. The method automates the computation of the singular and near singular precomputation integrals to a prescribed accuracy and determines the near/far interaction cutoff distance for each integration patch individually. The approach was used to solve scattering problems on PEC models of a sphere, toroid, and a complex CAD glider. The results from the three presented examples indicate that the proposed approach exhibits comparable performance regardless of whether the Gauss-Kronrod or Clenshaw-Curtis quadrature rule is employed.  The results confirm that the proposed approach delivers high-order accuracy and computational efficiency comparable to that of the fixed-grid CBIE method without requiring any geometry-specific calibration. Future work aims to address sharp-edge singularities to further extend applicability to broader classes of scattering problems.


\bibliographystyle{ieeetr}
\bibliography{3dadaptivebie-refs}

\end{document}